\documentclass[11pt,twoside,a4paper,reqno,fleqn]{amsart} 
\usepackage{amsmath,amssymb,amscd} 
\usepackage{times}             

\theoremstyle{plain} 
\newtheorem{theorem}{Theorem}[section] 
\newtheorem{proposition}[theorem]{Proposition} 
 
\newtheorem{lemma}[theorem]{Lemma} 
\newtheorem{corollary}[theorem]{Corollary}

\theoremstyle{definition} 
\newtheorem{definition}[theorem]{Definition} 
 
\theoremstyle{remark} 
\newtheorem{remark}[theorem]{Remark} 
\newtheorem{example}[theorem]{Example} 
 
\numberwithin{equation}{section} 
 
\newcommand{\N}{{\mathbb N}} 
 
\newcommand{\R}{{\mathbb R}} 
\newcommand{\Rp}{\R_+} 
\newcommand{\Rn}{\R^n} 
\newcommand{\C}{{\mathbb C}} 
 
\newcommand{\cB}{{\mathcal B}} 
 
\newcommand{\cL}{{\mathcal L}} 
 
\newcommand{\lDr}{\langle D_x\rangle} 
\newcommand{\lxr}{\langle\xi\rangle} 
\newcommand{\ls}{{l_\ast}} 
\newcommand{\bs}{{\beta_\ast}} 
\newcommand{\Op}{\operatorname{Op}} 
\newcommand{\cl}{{\textnormal{cl}}} 
\newcommand{\diag}{\operatorname{diag}}

\newcommand{\comp}{{\textnormal{comp}}} 
\newcommand{\bsen}{{\boldsymbol 1}_N} 
\newcommand{\bsenn}[1]{{\boldsymbol 1}_{{#1}}} 
\renewcommand{\Re}{\operatorname{Re}} 
\renewcommand{\Im}{\operatorname{Im}}

\newcommand{\AP}{{\em a~priori }\/}

\newcommand{\norm}[2]{\left\lVert{#1}\right\rVert_{#2}}   
\newcommand{\lp}[1]{L^{\!#1}}                         
\newcommand{\hn}[1]{H^{{#1}}}                         

\newcommand{\betrag}[1]{\langle{#1}\rangle} 

\newcommand{\comm}[2]{\left[{#1},{#2}\right]}  
\newcommand{\SP}[2]{\left({#1},{#2}\right)}    

\newcommand{\ve}{\varepsilon} 
\newcommand{\vr}{\varrho} 
\newcommand{\vp}{\varphi}

\renewcommand{\L}{{\cL}} 
\renewcommand{\S}{{\mathcal{S}}}

\newcommand{\del}{\partial}

\newcommand{\lap}{\bigtriangleup} 
 
\title[Weakly hyperbolic systems]{Energy estimates for  
  weakly \\ hyperbolic systems of the first order} 
\author{Michael Dreher} 
\address{Technische Universit\"at Bergakademie Freiberg,  
  Fakult\"at f\"ur Mathematik und Informatik, Institut f\"ur 
  Angewandte Analysis, D-09596 Freiberg, Germany} 
\email{dreher@math.tu-freiberg.de} 
\author{Ingo Witt}  
\address{Universit\"at Potsdam, Institut f\"ur Mathematik, Postfach 
  60\,15\,53, D-14415 Potsdam, Germany}   
\email{ingo@math.uni-potsdam.de} 
\subjclass{Primary: 35L80; Secondary: 35L30, 35S10}  
\keywords{Weakly hyperbolic systems, coefficients depending on the spatial
  variables, well-posedness of the Cauchy problem in Sobolev-type spaces, loss
  of regularity}
\date{\today} 
 
\begin{document} 

\begin{abstract} 
  For a class of weakly hyperbolic systems of the form $D_t-A(t,x,D_x)$, where
  $A(t,x,D_x)$ is a first-order pseudodifferential operator whose principal
  part degenerates like $t^\ls$ at time $t=0$, for some integer $\ls\geq1$,
  well-posedness of the Cauchy problem is proved in an adapted scale of
  Sobolev spaces. In addition, an upper bound for the loss of regularity that
  occurs when passing from the Cauchy data to the solutions is established. In
  examples, this upper bound turns out to be sharp.
\end{abstract} 
 
\renewcommand{\subjclassname}{\textup{2000} Mathematics Subject 
     Classification} 
\maketitle 
 
\setcounter{tocdepth}{1} 
\tableofcontents 
 
 
\section{Introduction} 
  

In this paper, we study the Cauchy problem for weakly hyperbolic
systems of the form
\begin{align} 
\label{syst} 
& 
\left\{\quad 
\begin{aligned} 
D_{t}U(t,x) &= A(t,x,D_{x})U(t,x)+F(t,x), 
\quad (t,x)\in (0,T)\times\R^{n}, 
\\ 
U(0,x) &=U_{0}(x),
\end{aligned} 
\right. 
\end{align} 
where $A(t,x,D_{x})$ is an $N\times N$ first-order pseudodifferential
operator. The precise assumptions on the symbol $A(t,x,\xi)$ are stated in
\eqref{FormA} below.
 
In order to motivate these assumptions, let us discuss an example. Systems of
the form \eqref{syst} arise, e.g., in converting $m$\/th-order partial
differential operators $P$ with principal symbol
\begin{equation}\label{PrincipSymb} 
\sigma^m(P)(t,x,\tau,\xi) = 
\prod_{h=1}^{m}(\tau-t^{l_{\ast}}\mu_{h}(t,x,\xi)), 
\quad \ls\geq1,
\end{equation} 
where $\mu_{h}\in C^{\infty}([0,T],S^{(1)})$ for $1\leq h\leq m$, into
first-order systems.

Assuming strict hyperbolicity for $t>0$, i.e., the $\mu_{h}$ are real-valued
and mutually distinct, it is well-known, see, e.g., {\sc Ivrii--Petkov}
\cite{IP74}, that the Cauchy problem for the operator $P$ is well-posed in
$C^{\infty}$ if and only if the lower-order terms satisfy so-called Levi
conditions. In case of~\eqref{PrincipSymb}, Levi conditions are expressed as
\begin{equation}\label{Levi}
P=\sum_{j+|\alpha|\leq m} 
a_{j\alpha}(t,x)\,t^{(j+(l_{\ast}+1)|\alpha|-m)^{+}}D_t^{j}D_x^{\alpha}, 
\end{equation}
with the coefficients $a_{j\alpha}(t,x)$ being smooth up to $t=0$.

Operators of the form \eqref{Levi} satisfying \eqref{PrincipSymb} are
particularly interesting because of two phenomena both occuring when passing
from the Cauchy data posed at $t=0$ to the solutions in the region $t>0$: One
is a certain loss of regularity and the other one is that the singularities
are possibly propagated in a non-standard fashion. These phenomena depend on
the lower-order terms of $P$ in a sensitive way.
 
One of the first examples in this direction was given by {\sc Qi} \cite{Qi},
\begin{equation}\label{QiEq} 
\left\{\quad 
\begin{aligned} 
& 
u_{tt}(t,x)-t^{2}u_{xx}(t,x)-(4k+1)u_{x}(t,x)=0, 
\quad (t,x)\in (0,T)\times\R, 
\\ 
& 
u(0,x)=\vp(x), \quad u_{t}(0,x)=0,   
\end{aligned} 
\right. 
\end{equation} 
where $k\in\N$. The solution to \eqref{QiEq} is 
\[ 
u(t,x)=\sum_{j=0}^{k}c_{jk}\,t^{2j}
\vp^{(j)}(x+t^{2}/2)
\] 
for certain coefficients $c_{jk}\neq0$. We see that $u(t,\cdot)$ for $t>0$ has
$k$ derivatives less compared to $\vp$. One actually loses $k$ derivatives for
any real number $k\geq-1/4$, as can be shown by an explicit representation of
the solution using special functions, see {\sc Taniguchi--Tozaki}
\cite{TaniguchiTozaki}. The parameter $k$ can even be a function $k(t,x)$ with
$k(0,x)\geq-1/4$ leading to a loss of regularity of $k(0,x)$, see {\sc
  Dreher}~\cite{DreherOsaka}. Further results concerning representation
formulas for the solutions and the propagation of singularities can be found
in {\sc Amano--Nakamura}~\cite{AmanoNakamura84}, {\sc
  Dreher--Reissig}~\cite{DreherReissig}, {\sc Hanges}~\cite{Han79}, {\sc
  Nakamura-Uryu}~\cite{NakamuraUryu}, {\sc Yagdjian}~\cite{YagdjianBuch}, {\sc
  Yoshikawa}~\cite{Yos77}. It is among the goals of this paper to establish
corresponding results for first-order weakly hyperbolic systems.
 
The first line of \eqref{QiEq} is converted into a first-order system by
setting
\[
U(t,x)= 
\begin{pmatrix} 
g(t,D_{x})u(t,x)\\ 
D_{t}u(t,x)   
\end{pmatrix}, 
\]
where 
\[
g(t,\xi) 
= 
\left( 
1-\chi 
\left( 
t^{2}\betrag{\xi}/2 
\right) 
\right) 
\betrag{\xi}^{1/2}
+
\chi 
\left( 
t^{2}\betrag{\xi}/2 
\right) 
t\betrag{\xi}, 
\]
and $\chi\in C^{\infty}(\overline{\R}_{+};\R)$ fulfills $\chi(t)=0$ if
$t\leq1/2$ and $\chi(t)=1$ if $t\geq1$. The symbol $g(t,\xi)$ will play an
important role later on. We then obtain
\begin{equation}\label{System1}  
D_{t}U(t,x)=A(t,x,D_{x})U(t,x), 
\quad  
(t,x)\in (0,T)\times\R, 
\end{equation}
where
\begin{equation}\label{grag}
A(t,x,\xi) 
= 
\chi 
\left( 
t^{2}\betrag{\xi}/2 
\right)
\left( 
t|\xi| 
\begin{pmatrix} 
0 & 1\\ 
1 & 0   
\end{pmatrix} 
- it^{-1} 
\begin{pmatrix} 
1 & 0 \\ 
b(t,x,\xi) & 0   
\end{pmatrix} 
\right) 
+A_{2}(t,x,\xi), 
\end{equation}
$b(t,x,\xi):=(4k(t,x)+1)\operatorname{sgn}\xi$, and $A_{2}(t,x,\xi)$ comprises
several terms of order zero, and other terms supported in the region
$t^{2}\lxr\leq 2$.

Generalizing \eqref{grag}, we are going to consider the Cauchy
problem~\eqref{syst} with operators $A(t,x,D_{x})$ whose symbols are of the
form
\begin{equation}\label{FormA}  
A(t,x,\xi) 
= 
\chi(\Lambda(t)\lxr) 
\left( 
\lambda(t)|\xi|A_{0}(t,x,\xi)-il_{\ast}t^{-1}A_{1}(t,x,\xi)   
\right) 
+ 
A_{2}(t,x,\xi), 
\end{equation} 
where $A_{0},\,A_{1}\in C^{\infty}([0,T],S^{(0)})$, $A_{2}\in
S^{-1,1;\lambda}+S^{0,0;\lambda}$ are $N\times N$ matrix-valued
pseudodifferential symbols, the function $\lambda(t)=t^\ls$ for the fixed
integer $\ls\geq1$ characterizes the \emph{kind of degeneracy\/} at $t=0$,
$\Lambda(t):=\int_{0}^{t}\lambda(t')\,dt'= \bs\,t^{\ls+1}$ is its primitive,
and $\bs:=1/(\ls+1)$. The symbol classes $S^{m,\eta;\lambda}$ for
$m,\,\eta\in\R$ that are closely related to the kind of degeneracy under
consideration will be introduced in Section~\ref{Sec21}.  In fact, these
symbol classes are characterized by two weight functions
\begin{align*}
   g(t,\xi)&:=(1-\chi(\Lambda(t)\lxr))\lxr^\bs + \chi(\Lambda(t)\lxr)
   \lambda(t)\lxr \\
   h(t,\xi)&:=(1-\chi(\Lambda(t)\lxr))\lxr^\bs + \chi(\Lambda(t)\lxr)t^{-1},
\end{align*}
where $m$ is the exponent of $g(t,\xi)$ and $\eta-m$ is the exponent of
$h(t,\xi)$. Note that $A(t,x,\xi)$ in \eqref{FormA} belongs to the class 
$S^{1,1;\lambda}$. 

We will present a symbolic calculus for matrices $A(t,x,\xi)$ of the form
\eqref{FormA} that for many purposes allows to argue on a purely algebraic
level, leading in this way to short and compact proofs.
 
We will also introduce function spaces $H^{s,\delta(x);\lambda}((0,T)\times
\R^n)$ to which the solutions $U(t,x)$ to \eqref{syst} belong. Here, $s\in\R$
is the Sobolev regularity with respect to $(t,x)$ for $t>0$, while $\delta=
\delta(x)$ is related to the loss of regularity at the point $x\in\R^{n}$. For
instance, for $s\in\N$, $\delta(x)=\delta$ being a constant, the space
$H^{s,\delta;\lambda}((0,T)\times\R^{n})$ consist of all functions $U(t,x)$
satisfying $k_{s-j,s+\delta}(t,x,D_x) D_t^jU\in L^{2}((0,T)\times\R^{n})$ for
$0\leq j\leq s$ and arbitrary $k_{m\eta}\in S^{m,m+\eta l_{\ast};\lambda}$.
The case of variable $\delta(x)$ will be discussed in detail in Section~3.
 
Our main result is the following:
 
\begin{theorem}\label{main} 
\textup{(a)} Assume the symbol $A(t,x,\xi)$ in~\eqref{FormA} is
symmetrizable-hyperbolic in the sense that there is a matrix $M_{0}\in
C^{\infty}([0,T],S^{(0)})$ such that $|\det M_{0}(t,x,\xi)|\geq c$ for some
$c>0$ and the matrix $(M_{0}A_{0}M_{0}^{-1})(t,x,\xi)$ is symmetric for all
$(t,x,\xi)\in [0,T] \times\Rn\times(\Rn\setminus0)$. 
Then there is a function $\delta\in \cB^\infty(\R^{n};\R)$ such that, for all
$s\geq0$, $U_{0}\in H^{s+\beta_{\ast}\delta(x)l_{\ast}}(\R^{n})$, and $F\in
H^{s,\delta(x);\lambda}((0,T)\times\R^{n})$, system \eqref{syst} possesses a
unique solution $U\in H^{s,\delta(x);\lambda}((0,T)\times\R^{n})$. Moreover,
we have the estimate
\begin{equation}\label{EE}
  \|U\|_{H^{s,\delta(x);\lambda}} \leq C\left(
  \|U_0\|_{H^{s+\bs\delta(x)\ls}}+\|F\|_{H^{s,\delta(x);\lambda}}\right)
\end{equation}
for a suitable constant $C=C(s,\delta,T)>0$. In particular, the loss of
regularity that is indeed independent of $s\geq0$ does not exceed
$\bs\delta(x)\ls$.
 
\textup{(b)} In\/ \textup{(a)}, we can choose any function
$\delta\in\cB^\infty(\R^{n};\R)$ for which there is matrix $M_1\in S^{(0)}$
such that the inequality
\begin{equation}\label{PRE}
\Re 
\left( 
M_{0}A_{1}M_{0}^{-1} 
+ 
\comm{M_{1}M_{0}^{-1}}{M_{0}A_{0}M_{0}^{-1}}   
\right)(0,x,\xi) 
\leq 
\delta(x)\bsen 
\end{equation}
holds for all $(x,\xi)\in\Rn\times(\Rn\setminus0)$. Here $[\;,\,]$ denotes the
commutator and $\Re Q:=(Q+Q^{\ast})/2$. 
\end{theorem} 

\begin{remark}\label{ergaen}
(a) Part (a) of Theorem~\ref{main} continues to hold if one solely assumes that
$A\in \Op S^{1,1;\lambda}$ and there is an invertible $M\in\Op
S^{0,0;\lambda}$ such that $\Im(MAM^{-1})\in\Op S^{0,1;\lambda}$,
cf.~Lemma~\ref{equiva}. In this situation, however, we have no simple formula
for $\delta\in\cB^\infty(\R^{n};\R)$.

(b) The loss of regularity for the weakly hyperbolic operator $P$ from
\eqref{Levi} equals $\bs(\delta(x)+m-1)\ls$, where
$\delta\in\cB^\infty(\R^{n};\R)$ is the function satisfying \eqref{PRE} for
the first-order systems arising by converting $P$.
\end{remark}
 
The paper is organized as follows: In Section 2, we introduce the symbol
classes $S^{m,\eta;\lambda}$ and certain subclasses
$\tilde{S}^{m,\eta;\lambda}$ thereof, where the latter contains symbols
$A(t,x,\xi)$ that possess "one and a half" principal symbols 
\[
  \sigma^{m}(A)\in t^{m(\ls+1)-\eta}\,C^\infty([0,T];S^{(m)}), 
  \quad \tilde{\sigma}^{m-1,\eta}(A)\in S^{(m-1)}. 
\]
A similar calculus, but differentiation with respect to $t$ is included in the
pseudodifferential action, was established by \textsc{Witt}~\cite{Wit03}.  In
case $\ls=1$, there is related work by {\sc Boutet de Monvel}~\cite{BdM74},
{\sc Joshi}~\cite{Jos98}, {\sc Yoshikawa}~\cite{Yos77}, and others.

Eq.~\eqref{FormA} actually defines the class $\tilde S^{1,1;\lambda}$, where
$\sigma^1(A)(t,x,\xi)=\lambda(t)|\xi|\, A_0(t,x,\xi)$,
$\tilde{\sigma}^{0,1}(A)(x,\xi)=-i\ls\,A_1(0,x,\xi)$ for $A(t,x,\xi)$ as given
there. According to Theorem~\ref{main}, $\sigma^1(A)(t,x,\xi)$,
$\tilde{\sigma}^{0,1}(A)(x,\xi)$ are exactly the symbols which are needed to
symmetrize system \eqref{syst} and to read off the loss of regularity,
respectively. In particular, for \textsc{Qi}'s example \eqref{System1}, we have
\[ 
\sigma^{1}(A)(t,x,\xi)= t|\xi|
\begin{pmatrix} 
0 & 1 \\ 
1 & 0   
\end{pmatrix}, 
\quad 
\tilde{\sigma}^{0,1}(A)(x,\xi) 
= 
-i 
\begin{pmatrix} 
1 & 0 \\ 
b(0,x,\xi) & 0   
\end{pmatrix}. 
\] 

We choose $M_0(t,x,\xi)=\frac{1}{2}\left(\begin{smallmatrix} 1 & -1 \\ 1
    & 1 \end{smallmatrix}\right)$, $(M_1M_0^{-1})(x,\xi)=\frac14
\left(\begin{smallmatrix} 0 & b-1 \\ b+1 & 0
  \end{smallmatrix}\right)(0,x,\xi)$ and obtain
\[
  \Re\left(
  M_{0}A_{1}M_{0}^{-1} + \comm{M_{1}M_{0}^{-1}}{M_{0}A_{0}M_{0}^{-1}}
  \right)(0,x,\xi)=\frac12 \begin{pmatrix} 1-\Re b & 0 \\ 
  0 & 1+\Re b
  \end{pmatrix}(0,x,\xi).
\]
This leads to a loss of regularity of $|\Re k(0,x)+\frac14|-\frac14$ for
Eq.~\eqref{QiEq}, see Remark~\ref{ergaen}~(b). Moreover, this result is
sharp. The reason that we decided to introduce $\delta\in\cB^\infty(\Rn;\R)$
in \eqref{PRE} via an estimate (rather than an equality) is that the factual
loss of regularity is Lipschitz as function of $x$, but may fail to be
of class $C^1$, as this example shows.

Section 3 is concerned with properties of the function spaces
$H^{s,\delta(x);\lambda}((0,T)\times\R^{n})$. We extend results of {\sc
  Dreher--Witt}~\cite{DW02} from the case of constant $\delta$ to the case of
variable $\delta(x)$.
Our main result Theorem~\ref{main} is then proved in Section 4. 
In Section 5, some special cases in which the \AP estimate \eqref{EE} is
employed are considered: differential systems, systems with characteristic
roots of constant multiplicity for $t>0$, and higher-order equations.
Choosing the matrix $M_{1}$ suitably, we will find that the upper bound for
the loss of regularity, as predicted by inequality \eqref{PRE}, coincides with
the actual loss of regularity, as known in special cases, see, e.g.,
\textsc{Nakamura--Uryu}~\cite{NakamuraUryu}. In a forthcoming paper, we will
provide lower bounds for the loss of regularity for system~\eqref{syst}, and
we will show that for a wide class of operators the \AP estimate given in the
present paper is sharp.

Finally, in an appendix we provide an estimate that is useful to bring the
remainder term $A_{2}(t,x,D_x)\in \Op S^{0,0;\lambda}+ \Op S^{-1,1;\lambda}
\subset \lp{\infty}((0,T),\Op S^\bs_{1,0})\cap t^{-1}\lp{\infty}((0,T),\Op
S^{0}_{1,0})$ under control.

 
\section{Symbol classes}\label{symbols} 
 
\subsection{The symbol classes \boldmath $S^{m,\eta;\lambda}$}\label{Sec21} 
 
 
In this section, we introduce the fundamental symbol classes $S^{m,\eta;
  \lambda}$ for $m,\,\eta\in\R$. For an $m$th-order symbol $a(t,x,\xi)$, the
belonging of $a$ to $S^{m,m;\lambda}$ expresses the fact that $\sigma^m(a)$
degenerates like $\lambda^m(t)$ at time $t=0$, and it expresses sharp Levi
conditions on the lower order terms as well. Note that corresponding symbol
estimates (involving the functions $\bar g$, $\bar h$ from \eqref{bar}) are
predicted by the definition of the function spaces $H^{s,\delta;\lambda}$ in
Section~\ref{spaces}. To be able to deal with operators that arise in reducing
\eqref{syst} with the help of the operator $\Theta$ from Lemma~\ref{Theta},
where the latter is zeroth-order for $t>0$, but of variable
order~$\bs\delta(x)\ls$ when restricted to time $t=0$, we further introduce
the symbol classes $S_+^{m,\eta;\lambda}$ as slightly enlarged versions of
$S^{m,\eta;\lambda}$, but for $m\in\R$, $\eta\in \cB^\infty(\Rn;\R)$.
 
In the sequel, all symbols $a(t,x,\xi)$ will take values in $N\times 
N$-matrices, for some $N\in\N$. 
 
Let 
\begin{equation}\label{bar} 
  \bar g(t,\xi):= \lambda(t)\lxr+\lxr^\bs, \quad  
  \bar h(t,\xi):= \bigl(t+\lxr^{-\bs}\bigr)^{-1}. 
\end{equation} 
 
\begin{definition} 
\textup{(a)} For $m,\,\eta\in\R$, the symbol class $S^{m,\eta;\lambda}$ 
consists of all $a\in C^\infty([0,T]\times\R^{2n};M_{N\times N}(\C))$ such 
that, for each $(j,\alpha,\beta)\in \N^{1+2n}$, there is a constant 
$C_{j\alpha\beta}>0$ with the property that 
\begin{equation}\label{symb_est} 
  \bigl|\partial_t^j\partial_x^\alpha\partial_\xi^\beta a(t,x,\xi)\bigr| 
  \leq C_{j\alpha\beta}\,\bar g(t,\xi)^m \bar h(t,\xi)^{\eta-m+j}  
  \lxr^{-|\beta|} 
\end{equation} 
for all $(t,x,\xi)\in [0,T]\times\R^{2n}$. 
 
\textup{(b)} For $m\in\R$, $\eta\in \cB^\infty(\Rn;\R)$, and $b\in\N$, the 
symbol class $S_{(b)}^{m,\eta;\lambda}$ consists of all $a\in 
C^\infty([0,T]\times\R^{2n};M_{N\times N}(\C))$ such that, for each 
$(j,\alpha,\beta)\in \N^{1+2n}$, there is a constant $C_{j\alpha\beta}>0$ with 
the property that 
\begin{equation}\label{symb_est2} 
  \bigl|\partial_t^j\partial_x^\alpha\partial_\xi^\beta a(t,x,\xi)\bigr| 
  \leq C_{j\alpha\beta}\,\bar g(t,\xi)^m \bar h(t,\xi)^{\eta(x)-m+j}  
  \bigl(1+|\log\bar h(t,\xi)|\bigr)^{b+|\alpha|}\lxr^{-|\beta|} 
\end{equation} 
for all $(t,x,\xi)\in [0,T]\times\R^{2n}$. Moreover, we set 
\[ 
  S_+^{m,\eta;\lambda}=\bigcup_{b\in\N} S_{(b)}^{m,\eta;\lambda}.  
\] 
\end{definition} 
 
As usual, we set  
\[ 
  S^{-\infty,\eta;\lambda} = \bigcap_{m\in\R}S^{m,\eta;\lambda}, 
\] 
see Proposition~\ref{properties} (a) and also (b). Similarly for 
$S_+^{-\infty,\eta;\lambda}$, $S_{(b)}^{-\infty,\eta;\lambda}$. 
 
\begin{remark} 
In view of $\bar g \bar h^\ls\sim \lxr$ and $\bar g \bar h^{-1} 
\sim 1+\Lambda(t)\lxr$, estimate \eqref{symb_est} is equivalent to 
\[ 
  \bigl|\partial_t^j\partial_x^\alpha\partial_\xi^\beta a(t,x,\xi)\bigr| 
  \leq C_{j\alpha\beta}'\,\bar g(t,\xi)^{m-|\beta|}  
  \bar h(t,\xi)^{\eta-m-|\beta|\ls+j}  
\] 
and  
\[ 
  \bigl|\partial_t^j\partial_x^\alpha\partial_\xi^\beta a(t,x,\xi)\bigr| 
  \leq C_{j\alpha\beta}''\,(1+\Lambda(t)\lxr)^m \bar h(t,\xi)^{\eta+j}  
  \lxr^{-|\beta|}, 
\] 
respectively. A similar remark applies to \eqref{symb_est2}. 
\end{remark} 
 
We discuss some examples of further use:
 
\begin{lemma}
Let $m,\,\eta\in\R$. Then\/\textup{:}
 
\textup{(a)} $\bar g^m\bar h^{\eta-m}\in S^{m,\eta;\lambda}$.  
 
\textup{(b)} For $a\in C^\infty([0,T];S^m)$ and $l\in\N$, $a\in 
S^{m,m(\ls+1)-l;\lambda}$ if and only if 
\[ 
  \partial_t^j a\big|_{t=0} \in S^{m-\bs(l-j)}, \quad 0\leq j\leq l-1. 
\] 
 
\textup{(c)} Let $\chi\in C^\infty(\overline{\R}_+;\R)$, $\chi(t)=0$ if 
$t\leq1/2$, $\chi(t)=1$ if $t\geq1$. Then  
\[ 
  \chi^+(t,\xi) := \chi(\Lambda(t)\lxr) \in S^{0,0;\lambda}, 
\] 
while $\chi^-(t,\xi):=1-\chi^+(t,\xi)\in S^{-\infty,0;\lambda}$. 
\end{lemma} 
 
In particular, from (a), (b) we infer 
\[ 
  \lambda(t)\lxr\in S^{1,1;\lambda}, \quad  
  (t+\lxr^{-\bs})^{-1}\in S^{0,1;\lambda}, \quad  \Lambda(t)\lxr\in  
  S^{1,0;\lambda}. 
\]  
 
In the next proposition, we list properties of the symbol classes 
$S^{m,\eta;\lambda}$ for $m,\,\eta\in\R$ (with proofs which are standard 
omitted): 
 
\begin{proposition}\label{properties} 
\textup{(a)} $S^{m,\eta;\lambda}\subseteq S^{m',\eta';\lambda}$ 
$\Longleftrightarrow$ $m\leq m'$, $\eta\leq \eta'$. 
 
\textup{(b)} Let $a\in S^{m,\eta;\lambda}$. Then $ \chi^+(t,\xi) a\in 
S^{m',\eta;\lambda}$ for some $m'<m$ implies $a\in S^{m',\eta;\lambda}$. In 
particular, if $a(t,x,\xi)=0$ for $\Lambda(t)\lxr\geq C$ and certain $C>0$, 
then $a\in S^{-\infty,\eta;\lambda}$. 
 
\textup{(c)} If $a\in S^{m,\eta;\lambda}$, then $\partial_t^j\partial_x^\alpha 
\partial_\xi^\beta a\in S^{m-|\beta|,\eta+j-|\beta|(\ls+1);\lambda}$. 
 
\textup{(d)} If $a\in S^{m,\eta;\lambda}$, $a'\in S^{m',\eta';\lambda}$, then 
$a\circ a'\in S^{m+m',\eta+\eta';\lambda}$ and 
\[ 
  a\circ a' = aa' \mod{S^{m+m'-1,\eta+\eta'-(\ls+1);\lambda}}, 
\] 
where $\circ$ denotes the Leibniz product with respect to $x$. 
 
\textup{(e)} If $a\in S^{m,\eta;\lambda}$, then  
$a^\ast\in S^{m,\eta;\lambda}$ and 
\[ 
  a^\ast(t,x,\xi) = a(t,x,\xi)^\ast \mod{S^{m-1,\eta-(\ls+1);\lambda}}, 
\] 
where $a^\ast$ is the\/ \textup{(}complete\/\textup{)} symbol of the formal
adjoint to $a(t,x,D_x)$ with respect to $L^2$.
 
\textup{(f)} If $a\in  S^{m,\eta;\lambda}([0,T]\times\R^{2n}; 
M_{N\times N}(\C))$ is elliptic in the sense that 
\[ 
  |\det a(t,x,\xi)|\geq c_1
  \left( 
  \bar{g}^m(t,\xi)\,\bar{h}^{\eta-m}(t,\xi) 
  \right)^N, \quad  
  (t,x,\xi)\in [0,T]\times \R^{2n},\,|\xi|\geq c_2 
\] 
for some $c_1,\,c_2>0$, then there is a symbol $a'\in S^{-m,-\eta;\lambda}$
with the property that
\[ 
  a\circ a' -1,\,a'\circ a-1 \in C^\infty([0,T];S^{-\infty}). 
\] 
Moreover, $a'$ is unique modulo $C^\infty([0,T];S^{-\infty})$ and
\[ 
  a' = a^{-1} \mod{S^{-m-1,-\eta-(\ls+1);\lambda}}. 
\] 
 
\textup{(g)} $\bigcap_{m,\eta} S^{m,\eta;\lambda}=C^\infty([0,T]; 
S^{-\infty})$. 
\end{proposition} 
 
Similar results hold for the classes $S_+^{m,\eta;\lambda}$ for $m\in\R$, 
$\eta\in \cB^\infty(\Rn;\R)$: 
 
\begin{proposition}\label{nnnm} 
\textup{(a)} $S_{(b)}^{m,\eta;\lambda}\subseteq S_{(b')}^{m',\eta'; \lambda}$ 
$\Longleftrightarrow$ $m\leq m'$, $\eta\leq \eta'$, and $b\leq b'$ if 
$\eta=\eta'$. 
 
\textup{(b)} $S^{m,\eta;\lambda}\subsetneq S_+^{m,\eta;\lambda} \subsetneq 
\bigcap_{\epsilon>0}S^{m,\eta+\epsilon;\lambda}$. 
 
\textup{(c)} If $a\in S_{(b)}^{m,\eta;\lambda}$, then 
$\partial_t^j\partial_x^\alpha \partial_\xi^\beta a\in 
S_{(b+|\alpha|)}^{m-|\beta|,\eta-|\beta|(\ls+1)+j;\lambda}$. 
  
\textup{(d)} If $a\in S_{(b)}^{m,\eta;\lambda}$, $a'\in 
S_{(b')}^{m',\eta';\lambda}$, then $a\circ a'\in 
S_{(b+b')}^{m+m',\eta+\eta';\lambda}$ and 
\[ 
  a\circ a' = aa' \mod{S_{(b+b'+1)}^{m+m'-1,\eta+\eta'-(\ls+1);\lambda}}. 
\] 
 
\textup{(e)} If $a\in S_{(b)}^{m,\eta;\lambda}$, then  
$a^\ast\in S_{(b)}^{m,\eta;\lambda}$ and 
\[ 
  a^\ast(t,x,\xi) = a(t,x,\xi)^\ast  
  \mod{S_{(b+1)}^{m-1,\eta-(\ls+1);\lambda}}. 
\] 
\textup{(f)} $S_{(0)}^{0,0;\lambda}\subset L^\infty((0,T);S_{1,\delta}^0)$ 
for any $0<\delta<1$.  
\end{proposition} 
 
From Proposition~\ref{nnnm}~(f) we conclude: 
 
\begin{corollary}\label{fer} 
$\Op S^{0,0;\lambda}\subset \Op S_{(0)}^{0,0;\lambda}\subset \cL(L^2)$.  
\end{corollary} 
 
 
\subsection{The symbol classes \boldmath $\tilde S^{m,\eta;\lambda}$} 
 
 
To establish precise upper bounds on the loss of regularity in
Theorem~\ref{main}~(b), we now refine the fundamental symbol classes
$S^{m,\eta;\lambda}$ to $\tilde S^{m,\eta;\lambda}$, where symbols
$a(t,x,\xi)$ in the latter class admit ``one and a half'' principal symbols
$\sigma^m(a)$, $\tilde\sigma^{m-1,\eta}(a)$. These principal symbols
enable us to read off the loss of regularity.
 
\begin{definition} 
For $m,\,\eta\in\R$, the class $\tilde S^{m,\eta;\lambda}$ consists of all  
$a\in S^{m,\eta;\lambda}$ that can be written in the form 
\begin{equation}\label{representation} 
  a(t,x,\xi) = \chi^+(t,\xi)\,t^{-\eta}\bigl( 
  a_0(t,x,t^{\ls+1}\xi) + a_1(t,x,t^{\ls+1}\xi)\bigr) + a_2(t,x,\xi), 
\end{equation} 
where 
\[ 
  a_0\in C^\infty([0,T];S^{(m)}), \quad a_1\in C^\infty([0,T];S^{(m-1)}), 
\] 
and $a_2\in S^{m-2,\eta;\lambda}+S^{m-1,\eta-1;\lambda}$. With 
$a(t,x,\xi)$ as in \eqref{representation} we associate the two symbols 
\begin{equation}\label{symbolST} 
  \sigma^m(a)(t,x,\xi) := t^{-\eta}\,a_0(t,x,t^{\ls+1}\xi), \quad 
  \tilde\sigma^{m-1,\eta}(a)(x,\xi) := a_1(0,x,\xi). 
\end{equation} 
\end{definition} 

\begin{remark} 
The symbol components $\chi^{+}(t,\xi)t^{-\eta}a_{j}(t,x,t^{l_{\ast}+1}\xi)$
in \eqref{representation} for $j=0,1$ belong to $S^{m-j,\eta;\lambda}$, while
$a_{2}(t,x,\xi)$ is regarded as remainder term.
\end{remark} 
 
For further use, we also introduce 
\[ 
  \tilde\sigma^{m,\eta}(a)(x,\xi) := a_0(0,x,\xi). 
\] 
Note that this symbol is directly derived from $\sigma^m(a)$. 
 
In the sequel, we shall employ the symbols 
\begin{align*} 
  g(t,\xi) &:=\chi^-(t,\xi)\,\lxr^{\bs}+\chi^+(t,\xi)\,\lambda(t)\lxr, \\ 
  h(t,\xi) &:=\chi^-(t,\xi)\,\lxr^{\bs}+\chi^+(t,\xi)\,t^{-1}. 
\end{align*} 
Note that $g\sim \bar g$, $h\sim \bar h$ so that the symbol estimates 
\eqref{symb_est} are not affected by this change.  
 
\begin{example}\label{exam2} 
(a) Let $m,\,\eta\in\R$. Then $g^m h^{\eta-m}\in \tilde 
S^{m,\eta;\lambda}$, 
\begin{equation}\label{frank}
  \sigma^m(g^m h^{\eta-m}) = t^{-\eta}\bigl(t^{\ls+1}|\xi|\bigr)^m, \quad 
  \tilde\sigma^{m-1,\eta}(g^m h^{\eta-m}) = 0. 
\end{equation}
 
(b) Let $a(t,x,\xi) := \sum_{|\alpha|\leq m}a_\alpha(t,x)\, 
t^{(|\alpha|(\ls+1)-m)^+}\xi^\alpha$, where $a_\alpha(t,x)\in 
\cB^\infty([0,T]\times\Rn)$ for $|\alpha|\leq m$. Then $a\in \tilde 
S^{m,m;\lambda}$, 
\begin{align*} 
  & \sigma^m(a) = \sum_{|\alpha|=m} 
    a_{j\alpha}(t,x)\,(t^\ls\xi)^\alpha, \\ 
  & \tilde\sigma^{m-1,m}(a) = \begin{cases} 
    \sum_{|\alpha|=m-1}a_{j\alpha}(0,x)\,\xi^\alpha & \text{if $m>1$}, \\ 
    0 & \text{if $m=0,\,1$.} 
  \end{cases} 
\end{align*} 
\end{example} 
 
The introduction of the principal symbols $\sigma^m(a)$, 
$\tilde\sigma^{m-1,\eta}(a)$ is justified by the next lemma: 
 
\begin{lemma} 
\textup{(a)} The symbols $\sigma^m(a)$, $\tilde\sigma^{m-1,\eta}(a)$ 
are well-defined. 
 
\textup{(b)} The short sequence 
\begin{equation}\label{short} 
\begin{CD} 
   0 @>>> S^{m-2,\eta;\lambda}+ S^{m-1,\eta-1;\lambda} @>>>  
   \tilde S^{m,\eta;\lambda} @>\left(\sigma^m,\tilde\sigma^{m-1,\eta}\right)>>
   \Sigma\tilde S^{m,\eta;\lambda} @>>> 0 
\end{CD} 
\end{equation} 
is exact, where $\Sigma \tilde S^{m,\eta;\lambda}:=\lambda^m(t)t^{-\eta+m}\, 
C^\infty([0,T];S^{(m)})\times S^{(m-1)}$ is the principal symbol space. 
\end{lemma} 
\begin{proof} 
For $a\in \tilde S^{m,\eta;\lambda}$ represented as in \eqref{representation}, 
we show that 
\[ 
  a\in S^{m-2,\eta;\lambda} + S^{m-1,\eta-1;\lambda} \quad  
  \Longleftrightarrow \quad a_0=0,\,a_1\big|_{t=0}=0. 
\] 
This gives (a) and also the exactness of the short sequence \eqref{short} in 
the middle. Since the surjectivity of the map 
$(\sigma^m,\tilde\sigma^{m-1,\eta})$ is obvious, the proof will then be 
finished. 
 
So, let us assume that $a_0\neq 0$ or $a_1\big|_{t=0}\neq 0$. If $a_0\neq0$,
then $|a|\geq C^{-1}\,g^m h^{\eta-m}$ for $\Lambda(t)\lxr\geq C$ in some conic
set, and $C>0$ sufficiently large. Hence, $a \notin S^{m-2,\eta;\lambda} +
S^{m-1,\eta-1;\lambda}$. If $a_0=0$, but $a_1\big|_{t=0}\neq 0$, then we write
\[ 
  a_1(t,x,\xi)= b_0(x,\xi) + tb_1(t,x,\xi),  
\] 
where $b_0\in S^{(m-1)}$, $b_1\in C^\infty([0,T],S^{(m-1)})$. But 
$\chi^+(t,\xi)\,t^{-\eta+1} b_1(t,x,t^{\ls+1}\xi)\in S^{m-1,\eta-1;\lambda}$, 
while $\chi^+(t,\xi)\, t^{-\eta} b_0(x,t^{\ls+1}\xi)\notin 
S^{m-2,\eta;\lambda} + S^{m-1,\eta-1;\lambda}$ in view of $b_0\neq0$. 
Hence, again, $a\notin S^{m-2,\eta;\lambda} + S^{m-1,\eta-1;\lambda}$.  
 
Now, assume $a_0= 0$ and $a_1\big|_{t=0}=0$. Write $a_1(t,x,\xi) = 
tb_1(t,x,\xi)$, where $b_1\in C^\infty([0,T],S^{(m-1)})$. Then 
\[ 
  a(t,x,\xi) = \chi^+(t,\xi)\,t^{-\eta+1} 
  b_1(t,x,t^{\ls+1}\xi) + a_2(t,x,\xi). 
\] 
But $\chi^+(t,\xi)\,t^{-\eta+1} b_1(t,x,t^{\ls+1}\xi)\in 
S^{m-1,\eta-1;\lambda}$, hence the claim. 
\end{proof}  
 
Finally, the next two results partially sharpen Proposition~\ref{properties}:
 
\begin{proposition}\label{krajt} 
\textup{(a)} If $a\in\tilde S^{m,\eta;\lambda}$, $a'\in \tilde 
S^{m',\eta';\lambda}$, then $a\circ a'\in \tilde S^{m+m',\eta+\eta'; 
\lambda}$ and 
\begin{align*} 
  \sigma^{m+m'}(a\circ a') &= \sigma^m(a)\,\sigma^{m'}(a'), \\ 
  \tilde\sigma^{m+m'-1,\eta+\eta'}(a\circ a') &=  
  \tilde\sigma^{m,\eta}(a)\,\tilde\sigma^{m'-1,\eta'}(a') + 
  \tilde\sigma^{m-1,\eta}(a)\,\tilde\sigma^{m',\eta'}(a'). 
\end{align*} 
 
\textup{(b)} If $a\in \tilde S^{m,\eta;\lambda}$, then $a^\ast\in 
\tilde S^{m,\eta;\lambda}$ and 
\[ 
  \sigma^m(a^\ast) = \sigma^m(a)^\ast, \quad  
  \tilde\sigma^{m-1,\eta}(a^\ast) = \tilde\sigma^{m-1,\eta}(a)^\ast. 
\] 
 
\textup{(c)} If $a\in \tilde S^{m,\eta;\lambda}([0,T]\times\R^{2n}; M_{N\times
  N}(\C))$ is elliptic, then $|\det \sigma^m(a)|\geq c\left(t^{(\ls+1)m-\eta}
  \, |\xi|^m\right){}^N$ for some $c>0$ and the symbol $a'$ from\/
\textup{Proposition~\ref{properties}~(f)} belongs to $\tilde
S^{-m,-\eta;\lambda}$. Moreover,
\[ 
  \sigma^{-m}(a')=\sigma^m(a)^{-1}, \quad  
  \tilde\sigma^{-m-1,-\eta}(a') = -\tilde\sigma^{m,\eta}(a)^{-1}\, 
  \tilde\sigma^{m-1,\eta}(a)\,\tilde\sigma^{m,\eta}(a)^{-1}. 
\] 
\end{proposition} 
\begin{proof} 
A straightforward computation. 
\end{proof} 
 
\begin{lemma}\label{fern} 
Let $a\in \tilde S^{m,\eta;\lambda}$ and $\eta=(\ls+1)m$. Then 
\[ 
  \partial_t a\in S^{m-1,\eta+1;\lambda}+S^{m,\eta;\lambda}.  
\] 
\end{lemma} 
\begin{proof} 
We have $\partial_t a\in \tilde S^{m,\eta+1;\lambda}$ and, in general, 
\[ 
  \tilde \sigma^{m,\eta+1;\lambda}(\partial_t a)=\left(m(\ls+1)-\eta\right)  
  \tilde \sigma^{m,\eta;\lambda}(a) 
\] 
Therefore, $\tilde \sigma^{m,\eta +1;\lambda}(\partial_t a)=0$ in case 
$\eta=(\ls+1)m$. The latter implies that $\partial_t a\in 
S^{m-1,\eta+1;\lambda}+S^{m,\eta;\lambda}$. 
\end{proof} 
 
\begin{remark} 
(a) For the reader's convenience, we summarize what vanishing of the single 
symbolic components for $a\in \tilde S^{m,\eta;\lambda}$ means: 
\begin{itemize} 
\item $\sigma^m(a)=0$, $\tilde\sigma^{m-1,\eta}(a)=0$ $\Longleftrightarrow$ 
  $a\in S^{m-2,\eta;\lambda} + S^{m-1,\eta-1;\lambda}$. 
\item $\sigma^m(a)=0$ $\Longleftrightarrow$ $a\in S^{m-1,\eta;\lambda}$. 
\item $\tilde\sigma^{m,\eta}(a)=0$ $\Longleftrightarrow$ $a\in 
  S^{m-1,\eta;\lambda} + S^{m,\eta-1;\lambda}$. 
\end{itemize} 
 
(b) Using the fact that asymptotic summation in the class $S^{m,\eta;\lambda}$ 
is possible one can introduce the class $S_\cl^{m,\eta;\lambda}$ of symbols 
$a\in S^{m,\eta;\lambda}$ obeying asymptotic expansions into double 
homogeneous components, and then it turns out that 
\[ 
  \tilde S^{m,\eta;\lambda} = S_\cl^{m,\eta;\lambda} + S^{m-2,\eta;\lambda}  
  + S^{m-1,\eta-1;\lambda}. 
\] 
The latter relation means that in $\tilde S^{m,\eta;\lambda}$ precisely the 
two symbolic components from \eqref{symbolST} survive. (Details on the class 
$S_\cl^{m,\eta;\lambda}$ will be published in a forthcoming paper 
\cite{DW04}.) 
\end{remark} 
 

\section{Function spaces}\label{spaces} 
 
In this section, we introduce the function spaces $H^{s,\delta;\lambda}$ for
$s\in\R$, $\delta\in \cB^\infty(\Rn;\R)$ and investigate their fundamental
properties. In case $s,\,\delta\in\R$, these function spaces were introduced
by \textsc{Dreher--Witt}~\cite{DW02} as abstract edge Sobolev spaces. Here, we
shall assume that the case of constant $\delta$ is known. Then the case of
variable $\delta$ is traced back to this previously known case. The key is the
invertibility of the operator $\Theta$, as stated in Lemma~\ref{Theta}.
 
\begin{definition}\label{func_sp} 
For $s\in\N$, $\delta\in \cB^\infty(\Rn;\R)$, the space $H^{s,\delta;\lambda}$
consists of all functions $u=u(t,x)$ on $(0,T)\times\Rn$ satisfying
\[ 
   (g^{s-j} h^{(s+\delta)\ls})(t,x,D_x)\,  
   D_t^j u \in L^2((0,T)\times\R^{n}), \quad 0\leq j\leq s.  
\]  
For general $s\in\R$, $\delta\in \cB^\infty(\Rn;\R)$, the space 
$H^{s,\delta;\lambda}$ is then defined by 
interpolation and duality.  
\end{definition} 

In particular, in case $s\geq0$, we have $(g^s h^{(s+\delta)\ls})(t,x,D_x)\, u
\in L^2((0,T)\times\R^{n})$ for any $u\in H^{s,\delta;\lambda}$.
 
\begin{remark}\label{later_use} 
(a) Strictly speaking, before Proposition~\ref{propertieST}~(a) we actually do 
not know that the spaces $H^{s,\delta;\lambda}$, firstly defined for 
$s\in\N$, interpolate. Therefore, it is only after 
Proposition~\ref{propertieST}~(a) that we get Lemma~\ref{Theta} and 
Proposition~\ref{map_prop} in their full strength. 
 
(b) Below we shall make use of Definition~\ref{func_sp} as follows:  
 
(i) For $s\in\N$, $\delta\in \cB^\infty(\Rn;\R)$, $u\in H^{s,\delta;\lambda}$ 
if and only if $g^{s-j}(t,D_x)D_t^ju\in H^{0,s+\delta;\lambda}$ for $0\leq 
j\leq s$. 
 
(ii) For $\delta\in \cB^\infty(\Rn;\R)$, $u\in H^{0,\delta;\lambda}$ if and 
only if $h^{\delta\ls}(t,x,D_x)u\in L^2$. 
\end{remark} 
 
For $K>0$, $\delta\in \cB^\infty(\Rn;\R)$, let $\lxr_K :=
(K^2+|\xi|^2)^{1/2}$, $\chi_K^+(t,\xi) := \chi(\Lambda(t)\lxr_K)$,
$\chi_K^-(t,\xi) := 1 - \chi_K^+(t,\xi)$, and
\[ 
  \Theta(t,x,\xi) = \Theta_{K,\delta}(t,x,\xi) := 
  \chi_K^-(t,\xi)\,\lxr_K^{\bs\delta(x)\ls} + 
  \chi_K^+(t,\xi)\,t^{-\delta(x)\ls}. 
\] 
Note that $\Theta(t,x,D_x)\in \Op S^{0,\delta(x)\ls;\lambda}_{(0)}$. 
 
\begin{lemma}\label{Theta} 
Given $\delta\in \cB^\infty(\Rn;\R)$, there is an $K_1>0$ such that the  
operator 
\begin{equation}\label{isom} 
  \Theta(t,x,D_x)\colon  
  H^{s,\delta';\lambda} \to H^{s,\delta'-\delta;\lambda}  
\end{equation} 
is invertible for all $s\in\R$, $\delta'\in \cB^\infty(\Rn;\R)$, and $K\geq 
K_1$. Moreover, $\Theta^{-1}\in \Op S^{0,-\delta(x)\ls;\lambda}_{(0)}$.  
\end{lemma} 
\begin{proof} 
Here, we will prove invertibility of the hypoelliptic operator 
$\Theta(t,x,D_x)$, for large $K>0$, and also the fact that 
$\Theta(t,x,D_x)^{-1}\in \Op S^{0,-\delta(x)\ls;\lambda}_{(0)}$. The proof 
is then completed with the help of the next proposition. 
 
The symbol $\Theta_{K,\delta}(t,x,\xi)$ belongs to the symbol class 
$S_+^{0,\delta(x)\ls;\lambda}$, but with parameter $K\geq K_0>0$. 
Similarly for $\Theta_{K,-\delta}(t,x,\xi)$. If $R_K':= \Theta_{K,\delta}\circ 
\Theta_{K,-\delta}-\Theta_{K,\delta}\Theta_{K,-\delta}$, then, for all 
$\alpha,\,\beta\in\N^n$ and certain constants $C_{\alpha\beta}>0$, 
\begin{multline*} 
  |\partial_x^\alpha \partial_\xi^\beta R_K'(t,x,\xi)|  
  \leq C_{\alpha\beta}\,\bigl(\lxr_K^\bs+\lambda(t)\lxr_K\bigr)^{-1} 
  (t+\lxr_K^{-\bs})^{-1} \\ 
  \times\, \bigl(1+ |\log(t+\lxr_K^{-\bs})|\bigr)^{1+|\alpha|} 
  \lxr_K^{-|\beta|}, \quad (t,x,\xi)\in[0,T]\times\R^{2n},\,K\geq K_0>0  
\end{multline*} 
(i.e., we have estimates \eqref{symb_est}, but with $\lxr$ replaced by
$\lxr_K$).  From the latter relation, it is seen that $R_K'(t,x,\xi)\to0$ in
$L^\infty((0,T);S^0)$ as $K\to\infty$, i.e., $R_K'(t,x,D_x)\to 0$ in
$\cL(L^2)$ as $K\to\infty$.
 
Now, let $R_K:=\Theta_{K,\delta} \circ \Theta_{K,-\delta}-1$, i.e.,
$R_K=R_K'+\Theta_{K,\delta}\Theta_{K,-\delta}-1$. Since
$(\Theta_{K,\delta}\Theta_{K,-\delta})(t,x,D_x)\to 1$ in $\cL(L^2)$ as
$K\to\infty$, it follows that $R_K(t,x,D_x)\to 0$ in $\cL(L^2)$ as
$K\to\infty$.  Thus, $\Theta_{K,-\delta}\circ (1+R_K)^{-1}$ is a right inverse
to $\Theta_{K,\delta}$, for large $K>0$. In a similar fashion, a left inverse
to $\Theta_{K,\delta}$ is constructed.
 
Moreover, $\Theta^{-1}=\Theta_{K,-\delta} \mod{\Op 
  S_+^{-\infty,-\delta(x)\ls-(\ls+1); \lambda}}$, as is seen from the 
constructions. 
\end{proof} 
 
\begin{proposition}\label{map_prop} 
For $m,\,s\in\R$, $\eta,\,\delta\in\cB^\infty(\Rn;\R)$, we have 
\begin{equation}\label{map} 
  \Op S^{m,\eta;\lambda}_{(0)} \subset \begin{cases} 
  \cL( H^{s,\delta;\lambda},H^{s-m,\delta+m+ 
  \frac{m-\eta}{\ls};\lambda}) & \textnormal{if $m\geq0$,} \\ 
  \cL( H^{s,\delta;\lambda},H^{s,\delta+ 
  \frac{m-\eta}{\ls};\lambda}) & \textnormal{if $m<0$.} \\ 
  \end{cases} 
\end{equation} 
\end{proposition} 
\begin{proof} 
We prove \eqref{map} in case $m\geq0$; the proof in case $m<0$ is similar. 
  
By interpolation and duality, we may assume that $s-m\in\N$. Then we have to 
show that, for $0\leq k\leq j\leq s-m$,  
\[ 
  h^{(s+\delta)\ls+m-\eta}g^{s-m-j}(D_t^{j-k}A)D_t^ku \in L^2   
\] 
provided $u\in H^{s,\delta;\lambda}$. We have 
\begin{equation}\label{help_inclusion} 
  h^{(s+\delta)\ls+m-\eta}g^{s-m-j}(D_t^{j-k}A)D_t^ku =  
  h^{m-\eta}g^{-m-j+k}(D_t^{j-k}A) 
  h^{(s+\delta)\ls}g^{s-k}D_t^k u + RD_t^ku
\end{equation} 
with $h^{m-\eta}g^{-m-j+k}(D_t^{j-k}A)\in\Op S^{-j+k,0;\lambda}_{(j-k)}$ and a
remainder $R\in \Op S_+^{s-j-1,(s-1)(\ls+1)+\delta\ls-k;\lambda}$. Now, $\Op
S^{-j+k,0;\lambda}_{(j-k)}\subset \Op S^{0,0;\lambda}_{(0)}$ and
$h^{(s+\delta)\ls}g^{s-k}D_t^k u\in L^{2}$ by assumption, i.e., the first
summand on the right-hand-side of \eqref{help_inclusion} belongs to $L^2$ by
virtue of Corollary~\ref{fer}. The second summand is rewritten as
\[ 
  RD_t^ku = Rg^{-s+k}(\Theta_{K,s+\delta})^{-1} 
  \Theta_{K,s+\delta}g^{s-k}D_t^k u 
\] 
for some large $K>0$, where $Rg^{-s+k}(\Theta_{K,s+\delta})^{-1}\in \Op 
S_+^{-j+k-1,-(\ls+1);\lambda}\subset \Op S^{0,0;\lambda}$ and again 
$\Theta_{K,s+\delta}g^{s-k}D_t^k u\in L^2$, i.e., also the 
second summand on the right-hand-side of \eqref{help_inclusion} belongs to 
$L^2$. 
\end{proof} 
 
In the following result, we summarize properties of the spaces 
$H^{s,\delta;\lambda}$. 
 
\begin{proposition}\label{propertieST} 
Let $s\in\R$, $\delta\in \cB^\infty(\Rn;\R)$. Then\/\textup{:} 
 
\textup{(a)} $\{H^{s,\delta;\lambda};\,s\in\R\}$ forms an interpolation scale 
of Hilbert spaces \textup{(}with the obvious Hilbert norms\/\textup{)} with 
respect to the complex interpolation method. 
 
\textup{(b)} $H^{s,\delta;\lambda}\big|_{\,(T',T)\times\Rn} =
H^s((T',T)\times\Rn)$ for all $0<T'<T$.
 
\textup{(c)} The space $C^{\infty}_{\comp}([0,T]\times\Rn)$ is 
dense in $H^{s,\delta;\lambda}$.  
 
\textup{(d)} For $s>1/2$, the map  
\begin{equation}\label{traces} 
  H^{s,\delta;\lambda} \to \prod_{j=0}^{[s-1/2]^{-}} H^{s+\bs\delta(x)\ls-
  \bs j-\bs/2}(\Rn), \quad u \mapsto \bigl(D_t^j u\big|_{t=0}
  \bigr)_{0\leq j\leq [s-1/2]^{-}},  
\end{equation} 
where $[s-1/2]^{-}$ is the largest integer strictly less than $s-1/2$, is
surjective.
 
\textup{(e)} $H^{s,\delta;\lambda}\subseteq H^{s',\delta';\lambda}$ if and 
only if $s\geq s'$, $s+\bs\delta\ls\geq s'+\bs\delta'\ls$. Moreover, the 
embedding $\{u\in H^{s,\delta;\lambda};\,\operatorname{supp}u\subseteq K\} 
\subseteq H^{s',\delta';\lambda}$ for some $K\Subset[0,T]\times\Rn$ 
is compact if and only if $s>s'$ and $s+\bs\delta(x)\ls> s'+\bs\delta'(x)\ls$ 
for all $x$ satisfying $(0,x)\in K$. 
\end{proposition} 
\begin{proof} 
For $s,\,\delta\in\R$, it is readily seen that Definition~\ref{func_sp}
coincides with that one given in \textsc{Dreher--Witt} \cite{DW02}. In this
case, proofs may be found there.  For variable $\delta=\delta(x)$, we
exemplarily verify (a), (d): To this end, we write
$H^{s,\delta;\lambda}=\Theta^{-1}H^{s,0;\lambda}$ for $s\in\R$, with $\Theta$
being the operator from Lemma~\ref{Theta}.
 
(a) Since $\{H^{s,0;\lambda};\,s\in\R\}$ is an interpolation scale, 
$\{H^{s,\delta;\lambda};\,s\in\R\}$ is also an interpolation scale with 
respect to the complex interpolation method. 
 
(d) Let $\gamma_j u:= D_t^j u\big|_{t=0}$. Then $\gamma_j\,\Theta u\in 
H^{s-\bs j -\bs/2}(\Rn)$ for $0\leq j\leq [s-1/2]^{-}$, since \eqref{traces} 
holds if $\delta=0$.  
 
Now, $H^{s,\delta;\lambda} \to \prod_{j=0}^{[s-1/2]^{-}}  
H^{s+\bs\delta(x)\ls-\bs j-\bs/2}(\Rn)$,  
$u\mapsto \bigl(\gamma_j u\bigr)_{0\leq j\leq [s-1/2]^{-}}$ follows 
from 
\[ 
  \gamma_j u = \bigl(\lDr_K^{\bs\delta(x)\ls}\bigr)^{-1}\gamma_j\Theta u, 
\] 
while the surjectivity of this map is implied by the reverse relation 
\[ 
  \gamma_j\,\Theta u = \lDr_K^{\bs\delta(x)\ls}\,\gamma_j u 
\] 
and the surjectivity of \eqref{traces} in case $\delta=0$. 
\end{proof} 
 
We also need the following results:
\begin{proposition}\label{asd} 
\textup{(a)} If $q(t,x,D_x)\in \Op S_{(0)}^{0,0;\lambda}$ is invertible
on $H^{s,\delta;\lambda}$ for some $s\in\R$, $\delta\in \cB^\infty(\Rn;\R)$,
then $q(t,x,D_x)$ is invertible on $H^{s,\delta;\lambda}$ for all $s\in\R$,
$\delta\in \cB^\infty(\Rn;\R)$ and
\[ 
  q(t,x,D_x)^{-1}\in \Op S_{(0)}^{0,0;\lambda}. 
\] 
 
\textup{(b)} Conversely, if $q_0\in C^\infty([0,T];S^{(0)})$ and $q_1\in 
S^{(-1)}$ are given, where $\left|\det q_0(t,x,\xi)\right|\geq c$ for all 
$(t,x,\xi)\in[0,T]\times\R^{2n}$ and a certain $c>0$, then there is an 
invertible operator $q(t,x,D_x)\in \Op \tilde S^{0,0;\lambda}$ in the sense  
of\/ \textup{(a)} such that 
\[ 
  \sigma^0(q)=q_0, \quad \tilde\sigma^{-1,0}(q)=q_1. 
\] 
\end{proposition} 
\begin{proof} 
(a) By conjugating the operator $q(t,x,D_x)$ with the inverse of $(g^s
h^{s\ls}\Theta)(t,x,D_x)$, where $\Theta(t,x,\xi)$ is as in Lemma~\ref{Theta},
we may suppose that $s=0$, $\delta=0$. From the invertibility of $q(t,x,D_x)$
on $H^{0,0;\lambda}$, we then conclude the ellipticity of the symbol
$q(t,x,\xi)$ in the standard fashion, i.e., we have $|\det q(t,x,\xi)|\geq c_1$
for all $(t,x,\xi)\in[0,T]\times\R^{2n}$, $|\xi|\geq c_2$, and some constants
$c_1,\,c_2>0$. By the analogue of Proposition~\ref{properties}~(f) for the
class $S_{(0)}^{0,0;\lambda}$, there is a symbol $q_1(t,x,\xi)\in
S_{(0)}^{0,0;\lambda}$ such that
\[ 
  q\circ q_1-1   \in C^\infty([0,T];S^{-\infty}). 
\] 
It follows that
\[
  q(t,x,D_x)\,q(t,x,D_x)^{-1} = q(t,x,D_x)\,q_1(t,x,D_x)  
  \mod C^\infty([0,T];\Op S^{-\infty}), 
\]
i.e., by multiplying both sides from the left by $q(t,x,D_x)^{-1}\in 
\Op S_{1,\bs}^0([0,T]\times\Rn\times\Rn)$,
\[ 
  q(t,x,D_x)^{-1} = q_1(t,x,D_x) \mod C^\infty([0,T];\Op S^{-\infty}) 
\]
and $q(t,x,D_x)^{-1}\in\Op S_{(0)}^{0,0;\lambda}$. 

(b) The rather long proof is deferred to Appendix A.2.
\end{proof} 
 
 
\section{Symmetrizable-hyperbolic systems} 
 
 
In this section we prove our main result Theorem~\ref{main}. 

 
\subsection{Reduction of the problem}  
 
 
For $A\in \Op\tilde S^{1,1;\lambda}$, throughout we shall adopt the notation 
\[ 
  \sigma^1(A)(t,x,\xi) = \lambda(t)|\xi|\,A_0(t,x,\xi), \quad  
  \tilde\sigma^{0,1}(A)(x,\xi) = -i\ls A_1(x,\xi), 
\] 
where $A_0\in C^\infty([0,T];S^{(0)})$, $A_1\in S^{(0)}$. Likewise, for the 
symmetrizer $M\in \Op\tilde S^{0,0;\lambda}$, we shall write  
\[ 
  \sigma^0(M)(t,x,\xi) = M_0(t,x,\xi), \quad \tilde\sigma^{-1,0}(M)(x,\xi) = 
  -i\ls|\xi|^{-1}M_1(x,\xi),  
\] 
where $M_0\in C^\infty([0,T];S^{(0)})$, $M_1\in S^{(0)}$.
Condition~\eqref{PRE} is
\begin{equation}\label{pre2}
  \Re\left(M_0A_1M_0^{-1} + \bigl[M_1M_0^{-1}, 
  M_0A_0M_0^{-1}\bigr]\right) \leq \delta(x)\bsen. 
\end{equation}  
 
\begin{remark} 
Because $M_0A_0M_0^{-1}$ is symmetric,  
\[ 
  \Re\bigl[M_1M_0^{-1},M_0A_0M_0^{-1}\bigr]= 
  i \bigl[\Im(M_1M_0^{-1}),M_0A_0M_0^{-1}\bigr], 
\] 
i.e., \eqref{pre2} amounts to choose $\Im(M_1M_0^{-1})$ appropriately.
\end{remark} 
 
\begin{lemma}\label{equiva} 
For system \eqref{syst} with $A(t,x,\xi)\in\tilde S^{1,1;\lambda}$, the
following conditions are equivalent\/\textup{:}
 
\textup{(a)} There is an $M_0\in C^\infty([0,T];S^{(0)})$ such that $|\det 
M_0(t,x,\xi)|\geq c$ for some $c>0$ and the matrix 
\[ 
  (M_0A_0M_0^{-1})(t,x,\xi)  
\] 
is symmetric for all $(t,x,\xi)\in [0,T]\times\Rn\times (\Rn\setminus0)$. 
 
\textup{(b)} There is an operator $M(t,x,D_x)\in\Op\tilde S^{0,0;\lambda}$ 
that is invertible on $L^2$ such that 
\[ 
  \Im (M A M^{-1}) \in \Op S^{0,1;\lambda}, 
\] 
i.e., $\Im \sigma^1(M A M^{-1}) = 0$. 
\end{lemma} 
\begin{proof} 
If (a) is fulfilled, let $M(t,x,D_x)\in\Op\tilde S^{0,0;\lambda}$ be
invertible such that $\sigma^0(M)(t,x,\xi)= M_0(t,x,\xi)$. Such an operator 
$M$ exists according to Proposition~\ref{asd} (b). Then we have that the
matrix 
\[
  \sigma^1(MAM^{-1})(t,x,\xi)=\lambda(t)|\xi|\,(M_0A_0M_0^{-1})(t,x,\xi)
\]
is symmetric for all $(t,x,\xi)$, i.e., $\sigma^1(\Im (M A M^{-1}))(t,x,\xi) =
0$ and $\Im (M A M^{-1}) \in \Op S^{0,1;\lambda}$.

Vice versa, if (b) is satisfied, then we can take $\sigma^0(M)(t,x,\xi)$ for
$M_0(t,x,\xi)$ in (a).
\end{proof} 
 
\begin{definition}\label{SymmHypDef} 
System \eqref{syst} is called symmetrizable-hyperbolic if the conditions of
Lemma~\ref{equiva} are fulfilled. It is called symmetric-hyperbolic if
$A_0(t,x,\xi)$ is already symmetric, i.e., $\Im A \in \Op S^{0,1;\lambda}$.
\end{definition} 
 
\begin{proposition}\label{reductProp} 
In the proof of\/ \textup{Theorem~\ref{main}}, we can assume that  
\begin{equation}\label{reduc} 
  A(t,x,\xi) = \chi^+(t,\xi)\left(\lambda(t)|\xi|\,A_0(t,x,\xi) - i\ls  
  t^{-1} A_1(x,\xi)\right) + A_2(t,x,\xi),  
\end{equation} 
where $A_0\in C^\infty([0,T];S^{(0)})$, $A_0 = A_0^\ast$, $A_1\in S^{(0)}$, 
\[ 
  \Re A_1(x,\xi) \leq 0, 
\] 
and $A_2\in S^{-1,1;\lambda} + S^{0,0;\lambda}$\textup{;} and $\delta=0$. 
\end{proposition} 
\begin{proof} 
Note that \eqref{reduc} means $\sigma^1(A)(t,x,\xi)$ is symmetric, while $\Im 
\tilde\sigma^{0,1}(A)(x,\xi)\geq0$.  
 
Let the assumptions of Theorem~\ref{main} be satisfied. In particular, let 
$\delta\in\cB^\infty(\Rn;\R)$ satisfy \eqref{PRE}. We reduce \eqref{syst} in 
two steps. 
 
(a) Using the symmetrizer $M\in\Op \tilde S^{0,0;\lambda}$, that is an 
isomorphism from $H^{s,\delta;\lambda}$ onto $H^{s,\delta;\lambda}$ for all 
$s\in\R$ by Proposition~\ref{map_prop}, while $M(0,x,D_x)$ is an isomorphism 
from $H^s(\Rn)$ onto $H^s(\Rn)$, instead of \eqref{syst} we consider the 
equivalent system satisfied by $V:=MU$: 
\begin{equation}\label{syst2} 
\left\{ \enspace 
\begin{aligned} 
  D_t V(t,x) &= B(t,x,D_x) V(t,x) + G(t,x), \quad (t,x)\in (0,T)\times\Rn, \\ 
  V(0,x) &= V_0(x). 
\end{aligned}, 
\right. 
\end{equation} 
where $B=MAM^{-1}+(D_tM)M^{-1}$, $V_0 = M(0,x,D_x)U_0$, $G=MF$. 
 
We have $B\in\Op\tilde S^{1,1;\lambda}$, $\sigma^1(B)=\lambda(t)|\xi|\, 
(M_0A_0M_0^{-1})(t,x,\xi)$,  
\[ 
  \tilde\sigma^{0,1}(B)=\tilde\sigma^{0,1}(MAM^{-1})= 
  -i\ls\left(M_0A_1M_0^{-1} + \bigl[M_1M_0^{-1},M_0A_0M_0^{-1}\bigr]\right) 
\] 
according to the composition rules in Proposition~\ref{krajt}. In the last 
line, it was employed that $(D_t M)M^{-1}\in\Op\tilde S^{0,1;\lambda}$, 
$\tilde\sigma^{0,1}((D_t M)M^{-1})=0$ by virtue of Lemma~\ref{fern}. 
 
Thus, we can assume that $A_0(t,x,\xi)$ is symmetric, $M_0(t,x,\xi)=\bsen$, 
and $M_1(x,\xi)=0$ in Theorem~\ref{main}. In this first reduction, $\delta$ 
has not been changed. 
 
(b) Now assume $A_0(t,x,\xi)$ is symmetric, $M_0(t,x,\xi)=\bsen$, and
$M_1(x,\xi)=0$. Then using the operator $\Theta$ from Lemma~\ref{Theta}, that
is an isomorphism from $H^{s,\delta;\lambda}$ onto $H^{s,0;\lambda}$ for all
$s\in\R$, while $\Theta(0,x,D_x)=\langle D_x\rangle_{K}^{\bs\delta(x)\ls}$ is
an isomorphism from $H^{s+\bs\delta(x)\ls}(\Rn)$ onto $H^s(\Rn)$, instead of
\eqref{syst} we consider the equivalent system satisfied by $V:=\Theta U$:
\begin{equation}\label{syst3} 
\left\{ \enspace 
\begin{aligned} 
  D_t V(t,x) &= B(t,x,D_x) V(t,x) + G(t,x), \quad (t,x)\in (0,T)\times\Rn, \\ 
  V(0,x) &= V_0(x). 
\end{aligned}, 
\right. 
\end{equation} 
where this time $B=\Theta A\Theta^{-1}+(D_t\Theta)\Theta^{-1}$, $V_0 = 
\Theta(0,x,D_x) U_0$, $G=\Theta F$. 
 
By Lemma~\ref{DD} below, $\Theta A \Theta^{-1} + (D_t\Theta)\Theta^{-1}\in 
\Op\tilde S^{1,1;\lambda}$, $\sigma^1(\Theta A \Theta^{-1} + 
(D_t\Theta)\Theta^{-1}) = \lambda(t)|\xi|A_0$, 
\[ 
  \tilde\sigma^{0,1}(\Theta A \Theta^{-1} + (D_t\Theta)\Theta^{-1})(x,\xi) 
  = -i\ls\left(A_1(x,\xi) - \delta(x)\right)\bsen. 
\] 
Thus we can, in addition, assume that $\Re A_1\leq0$. This second reduction 
changes $\delta$ to zero. 
\end{proof}

\begin{lemma}\label{DD} 
Let $\Theta$ be as in\/ \textup{Lemma~\ref{Theta}}. Then $\Theta A  
\Theta^{-1} + (D_t\Theta)\Theta^{-1}\in \Op\tilde S^{1,1;\lambda}$ and 
\begin{align*} 
  \sigma^1(\Theta A \Theta^{-1} + (D_t\Theta)\Theta^{-1}) 
  &= \sigma^1(A), \\ 
  \tilde\sigma^{0,1}(\Theta A \Theta^{-1} + (D_t\Theta)\Theta^{-1}) 
  &= \tilde\sigma^{0,1}(A) + i\delta(x)\ls\bsen. 
\end{align*} 
\end{lemma} 
\begin{proof}  
We have $\Theta A \Theta^{-1}\in \Op\tilde S^{1,1;\lambda}$ and  
$\sigma^1(\Theta A \Theta^{-1})=\sigma^1(A)$, $\tilde\sigma^{0,1}(\Theta A 
\Theta^{-1})=\tilde\sigma^{0,1}(A)$ because of 
\[ 
  \Theta\circ A\circ \Theta^{-1} = \Theta_{K,\delta} A \Theta_{K,-\delta} 
  = A \mod S_+^{-\infty,-\ls;\lambda}. 
\] 
Furthermore,  
\[ 
  (D_t\Theta)\circ \Theta^{-1} = \left(D_t\Theta_{K,\delta}\right) 
  \Theta_{K,-\delta} \mod{S_+^{-1,-\ls;\lambda}\subset S^{-1,1;\lambda}}; 
\] 
so we consider the product 
$\left(D_t\Theta_{K,\delta}\right) \Theta_{K,-\delta}$: 
\begin{align*} 
  \left(D_t\Theta_{K,\delta}\right) 
  \Theta_{K,-\delta} &= \Bigl(\lambda(t)\lxr_K\, 
  \chi'\left(\Lambda(t)\lxr_K\right) 
  \bigl(t^{-\delta(x)\ls} -\lxr_K^{\bs\delta(x)\ls}\bigr) \\ 
  & \qquad -\delta(x)\ls\,\chi_K^+(t,\xi)t^{-\delta(x)\ls-1}\Bigr)\, 
  \Bigl(\chi_K^-(t,\xi)\,\lxr_K^{-\bs\delta(x)\ls} +  
  \chi_K^+(t,\xi)\,t^{\delta(x)\ls}\Bigr) \\ 
  &= \lambda(t)\lxr_K\, 
  \chi'\left(\Lambda(t)\lxr_K\right)\chi_K^-(t,\xi) 
  \bigl((c_1\Lambda(t)\lxr_K)^{-\bs\delta(x)\ls} - 1\bigr) \\ 
  &\qquad + \lambda(t)\lxr_K\, 
  \chi'\left(\Lambda(t)\lxr_K\right)\chi_K^+(t,\xi) 
  \bigl(1-(c_1\Lambda(t)\lxr_K)^{\bs\delta(x)\ls}\bigr) \\ 
  &\qquad -\delta(x)\ls\,\chi_K^+(t,\xi)\chi_K^-(t,\xi) 
  \lambda(t)\lxr_K(c_1\Lambda(t)\lxr_K)^{-\bs\delta(x)\ls-1} \\ 
  &\qquad -\delta(x)\ls\left(\chi_K^+(t,\xi)\right)^2t^{-1} 
\end{align*} 
with $c_1=\ls+1$. The first three summands on the right-hand side belong to
$S^{-\infty,1;\lambda}$, since we have, e.g., $\chi'(t)\left(1-
  \chi(t)\right)\in C_\comp^\infty(\Rp)$; thus, $d_1\leq \Lambda(t)\lxr_K\leq
d_2$ for certain constants $0<d_1<d_2$ on the support of the first summand and
the derivatives of $(c_1\Lambda(t)\lxr_K)^{-\bs\delta(x)\ls}$ with respect to
$x$ do not produce logarithmic terms in the estimates.
 
Thus, we obtain 
\[ 
  (D_t\Theta)\circ \Theta^{-1} = i\delta(x)\ls\, 
  \chi^+(t,\xi)\,t^{-1} \mod{S^{-1,1;\lambda}}, 
\] 
i.e., $(D_t\Theta)\Theta^{-1}\in\Op\tilde S^{1,1;\lambda}$ and 
$\sigma^1((D_t\Theta)\Theta^{-1})=0$, $\tilde\sigma^{0,1}((D_t\Theta) 
\Theta^{-1})=i\delta(x)\ls\bsen$, as required. 
\end{proof} 
 
 
\subsection{Proof of Theorem~\ref{main}}\label{sec42} 


We now come to the proof of the main theorem. We divide this proof into
three steps. Thereby, we always assume the reductions made in
\textup{Proposition~\ref{reductProp}}.

\textbf{First step (Basic \AP estimate).}\enspace
Each solution $U$ to system~\eqref{syst} satisfies the \AP estimate \eqref{EE}
in case $s=0$, i.e.,
\begin{equation}\label{APs0Prop}
\norm{U}{H^{0,0;\lambda}((0,T)\times\R^{n})} 
\leq 
C 
\left( 
\norm{U_{0}}{L^2(\R^{n})} 
+ 
\norm{F}{H^{0,0;\lambda}((0,T)\times\R^{n})}   
\right),
\end{equation}
where $C=C(T)>0$.

\begin{proof} 
First recall that $H^{0,0;\lambda}((0,T)\times\R^{n})=L^2((0,T)\times\Rn)$.

Rewrite~\eqref{syst} in the form $(\del_{t}-B)U=iF$, where
\begin{align*}
 &B(t,x,\xi)=iA(t,x,\xi)=B_{1}(t,x,\xi)+B_{r}(t,x,\xi), \\
 &B_{1}(t,x,\xi) - i\chi^+(t,\xi)\left(\lambda(t)|\xi|\,A_0(t,x,\xi) - i\ls  
  t^{-1} A_1(x,\xi)\right)\in S^{-1,1;\lambda},
\end{align*}
and $B_{r}\in S^{0,0;\lambda}$. By construction,
\[ 
(B_{1}+B_{1}^{\ast})(t,x,\xi)\leq 2q(t,\xi)\bsen, 
\] 
where $q(t,\xi)=Cg(t,\xi)^{-1}h(t,\xi)^2$ and $\int_{0}^{t}q(t',\xi)\,dt'\in
\lp{\infty}((0,T),S^{0}_{1,0})$. From Lemma~\ref{UsefulLem}, we infer
\begin{equation}\label{nmn} 
\norm{U(t,\cdot)}{\lp{2}(\R^{n})}^2 
\leq 
C 
\left( 
\norm{U_{0}}{\lp{2}(\R^{n})}^2 
+ 
\int_{0}^{t} 
\norm{F(t',\cdot)}{\lp{2}(\R^{n})}^2 
\,dt'
\right). 
\end{equation}
Integrating this inequality over the time interval $(0,T)$ yields the desired
estimate \eqref{APs0Prop}.
\end{proof} 


\textbf{Second step (\AP estimate of higher-order derivatives).}\enspace
Each solution $U$ to system~\eqref{syst} satisfies the \AP estimate \eqref{EE}
in case $s>0$, i.e.,
\begin{equation}\label{APsProp}
\norm{U}{H^{s,0;\lambda}((0,T)\times\R^{n})} 
\leq 
C 
\left( 
\norm{U_{0}}{H^s(\R^{n})} 
+ 
\norm{F}{H^{s,0;\lambda}((0,T)\times\R^{n})}   
\right),
\end{equation}
where $C=C(s,T)>0$.

\begin{proof}
For any $s\in\R$, we have $U\in H^{s+1,0;\lambda}$ if and only if
$g(t,D_x)h^\ls(t, D_x)U,\,h^\ls(t,D_x)D_tU\in H^{s,0;\lambda}$. Moreover, the
vector $\left(g(t,D_x) h^\ls(t,D_x)U,\, h^\ls(t,D_x) D_tU\right)^T$ is a
solution to the Cauchy problem
\begin{equation}  
\left\{\quad 
\begin{aligned} 
  D_t\binom{gh^\ls U}{h^\ls D_tU} & = \begin{pmatrix} 
   A_{00} & 0 \\ A_{10} & A_{11} \end{pmatrix} \binom{gh^\ls U}{h^\ls D_tU}  
   + \binom{gh^\ls F}{D_t h^\ls F}, \\ 
  \binom{gh^\ls U}{h^\ls D_tU}(0,x) &= \binom{\lDr\, U_0(x)}{ 
   \lDr^{1-\bs}\left(A(0,x,D_{x})U_0(x)+F(0,x)\right)}, 
\end{aligned} 
\right. 
\end{equation} 
where 
\begin{align*} 
  A_{00} &= gh^\ls A (gh^\ls)^{-1} + (D_tg)g^{-1} + \ls (D_th)h^{-1},\\ 
  A_{10} &= \bigl[h^\ls(D_tA)+\ls (D_th)h^{\ls-1}A\bigr](gh^\ls)^{-1},\\ 
  A_{11} &= h^\ls A h^{-\ls}. 
\end{align*} 
By Lemma~\ref{fg0} below, induction on $s\in\N$, and interpolation in
$s\geq0$, we then deduce the second step from the first one.
\end{proof}
 
\begin{lemma}\label{fg0} 
We have $\left(\begin{smallmatrix} A_{00} & 0 \\ A_{10} & A_{11} 
\end{smallmatrix}\right)\in \Op\tilde S^{1,1;\lambda}$ and 
\begin{align}   
 \sigma^1\left(\begin{pmatrix} A_{00} & 0 \\ A_{10} & A_{11} 
  \end{pmatrix}\right) &=  \begin{pmatrix} \sigma^1(A) &  
  0 \\ 0 & \sigma^1(A)\end{pmatrix}, \label{fg10} \\
 \tilde\sigma^{0,1}\left(\begin{pmatrix} A_{00} & 0 \\ A_{10} &  
  A_{11} \end{pmatrix}\right) &=  \begin{pmatrix}  
   \tilde\sigma^{0,1}(A) & 0 \\ 0 & \tilde\sigma^{0,1}(A)\end{pmatrix}. 
   \label{fg11}
\end{align} 
In particular, $\left(\begin{smallmatrix} A_{00} & 0 \\ A_{10} & A_{11} 
\end{smallmatrix}\right)$ fulfills the same assumptions as $A\in  \Op\tilde 
S^{1,1;\lambda}$ does, but for $(2N)\times(2N)$ matrices. 
Furthermore,
\begin{equation}\label{fg2} 
 \binom{\lDr\, U_0}{\lDr^{1-\bs}\left(A(0)U_0+F(0)\right)}\in H^s(\Rn), \quad 
 \binom{gh^\ls F}{D_t h^\ls F}\in H^{s,0;\lambda} 
\end{equation} 
provided that $U_0\in H^{s+1}(\Rn)$, $F\in H^{s+1,0;\lambda}$. 
\end{lemma} 
\begin{proof} 
A straightforward calculation using Proposition~\ref{krajt} and \eqref{frank} 
gives $gh^\ls A (gh^\ls)^{-1}\in \Op\tilde S^{1,1;\lambda}$, 
\[ 
  \sigma^1(gh^\ls A (gh^\ls)^{-1})=\sigma^1(A), \quad  
  \tilde\sigma^{0,1}(gh^\ls A (gh^\ls)^{-1})=\tilde\sigma^{0,1}(A), 
\] 
$(D_tg)g^{-1},\,(D_th)h^{-1}\in \Op\tilde S^{0,1;\lambda}$, 
\[ 
  \tilde\sigma^{0,1}((D_tg)g^{-1})=-i\ls, \quad  
  \tilde\sigma^{0,1}((D_th)h^{-1})=i, 
\] 
$h^\ls(D_tA)(gh^\ls)^{-1},\,(D_th)h^{\ls-1}A(gh^\ls)^{-1}\in\Op\tilde 
S^{0,1;\lambda}$ 
\[ 
  \tilde\sigma^{0,1}(h^\ls(D_tA)(gh^\ls)^{-1})=-i\ls|\xi|^{-1} 
  \tilde\sigma^{1,1}(A), \quad \tilde\sigma^{0,1}((D_th) 
  h^{\ls-1}A(gh^\ls)^{-1})=i|\xi|^{-1}\sigma^{1,1}(A), 
\] 
and $h^\ls A h^{-\ls}\in \Op\tilde S^{1,1;\lambda}$, 
\[ 
  \sigma^1(h^\ls A h^{-\ls})=\sigma^1(A), \quad  
  \tilde\sigma^{0,1}(h^\ls A h^{-\ls})=\tilde\sigma^{0,1}(A). 
\] 
Thus, \eqref{fg10}, \eqref{fg11} hold. Moreover, \eqref{fg2} is obvious. 
\end{proof} 
  

\textbf{Third step (Existence and uniqueness).}\enspace
For all $U_{0}\in \hn{s}(\R^{n})$, $F\in H^{s,0;\lambda}((0,T)\times\R^{n})$,
where $s\geq0$, system~\eqref{syst} possesses a unique solution $U\in
H^{s,0;\lambda}((0,T)\times\R^{n})$ satisfying the \AP estimate
\eqref{APsProp}.
 
\begin{proof} 
Let $s\geq 1$, the general case then follows by density arguments. By 
Proposition~\ref{propertieST}~(c), we may suppose that  
$U_{0}\in C^{\infty}_{\comp}(\R^{n})$, 
$F\in C^{\infty}_{\comp}([0,T]\times\R^{n})$. 
 
We replace the operator $A(t,x,D_{x})$ by $A_{\ve}(t,x,D_{x})$ for $0<\ve\leq
1$, where
\begin{align*}
A_{\ve}(t,x,\xi)
&=
\chi^{+}(t,\xi)
\left(
\lambda(t)|\xi|A_{0}(t,x,\xi)
-
il_{\ast}(t+\ve)^{-1}A_{1}(x,\xi)
\right)
+
A_{2\ve}(t,x,\xi), \\
A_{2\ve}(t,x,\xi)
&=
\frac{t+\lxr^{-\bs}}{t+\lxr^{-\bs}+\ve}\,A_{2}(t,x,\xi).
\end{align*}
The system $D_{t}-A_{\ve}(t,x,D_x)$ is symmetrizable-hyperbolic with the
lower-order term belonging to the space $\lp{\infty}((0,T),S^{0}_{1,0})$.
Therefore, the Cauchy problem
\begin{align*} 
& 
\left\{\quad 
\begin{aligned} 
D_{t}U_{\ve}(t,x) &= A_{\ve}(t,x,D_{x})U_{\ve}(t,x)+F(t,x), 
\quad (t,x)\in (0,T)\times\R^{n}, 
\\ 
U_{\ve}(0,x) &=U_{0}(x), 
\end{aligned} 
\right. 
\end{align*} 
possesses a unique solution $U_{\ve}\in
C^{\infty}([0,T],\hn{\infty}(\R^{n}))$, see {\sc Taylor}~\cite{Tay81}.
 
The set $\{A_{\ve}\colon 0<\ve\leq1\}$ is bounded in
$\tilde{S}^{1,1;\lambda}$. Hence, the second step provides an estimate
\[ 
\norm{U_{\ve}}{H^{s,0;\lambda}((0,T)\times\R^{n})} 
\leq 
C 
\left( 
\norm{U_{0}}{\hn{s}(\R^{n})} 
+ 
\norm{F}{H^{s,0;\lambda}((0,T)\times\R^{n})}   
\right). 
\] 
that holds uniformly in $0<\ve\leq1$. Furthermore, the set
$\bigl\{(A_{\ve}-A_{\ve'})/(\ve-\ve')\colon 0<\ve'<\ve\leq 1\bigr\}$ is
bounded in $S^{0,2;\lambda}$. From the first step as well as
Propositions~\ref{map_prop}, \ref{propertieST}~(e), we deduce that
\begin{align*} 
\norm{U_{\ve}-U_{\ve'}}{H^{0,0;\lambda}} 
& 
\leq 
C\norm{(A_{\ve}-A_{\ve'})U_{\ve}}{H^{0,0;\lambda}} 
\\ 
& 
\leq 
C(\ve-\ve')\norm{U_{\ve}}{H^{0,2/l_{\ast};\lambda}} 
\leq 
C(\ve-\ve')\norm{U_{\ve}}{H^{1,0;\lambda}}
\end{align*}
for $0<\ve'<\ve\leq 1$. Since $s\geq 1$, this implies that $U_{\ve}$
converges to some limit $U$ in the space $H^{0,0;\lambda}$ as $\ve\to+0$. 

The rest of the proof is standard.
\end{proof} 
 
 
\section{Applications}\label{exam} 
 
We discuss three examples demonstrating the value of Theorem~\ref{main}.

\subsection{Differential systems}  
 
 
Differential systems of the form \eqref{syst} with $A(t,x,\xi)$
from~\eqref{FormA} are of restricted interest, because a lower-order term as
described by the term $\chi^+(t,\xi)t^{-1}A_{1}(t,x,\xi)$ cannot occur.
Hence, the loss of regularity is always zero.
 
Consider the operator 
\begin{equation}\label{diff_system} 
  L = D_t + \sum_{j=1}^n t^\ls a_j(t,x) D_{x_j} + a_0(t,x), 
\end{equation} 
where $a_j\in \cB^\infty([0,T]\times\Rn;M_{N\times N}(\C))$ for $0\leq j\leq
n$. With $A(t,x,\xi) := - \sum_{j=1}^n t^\ls a_j(t,x)\xi_j-a_0(t,x)$, 
\[ 
  \sigma^1(A)(t,x,\xi) = -\lambda(t)|\xi|\,\sum_{j=1}^n a_j(t,x) 
  \frac{\xi_j}{|\xi|}, \quad \tilde\sigma^{0,1}(A)(x,\xi) = 0. 
\] 

\begin{proposition} 
Let the differential system \eqref{diff_system} be symmetrizable-hyperbolic. 
Then, for all $s\geq0$, $U_0\in H^s(\Rn)$, and  
$F\in H^{s,0;\lambda}((0,T)\times\R^{n})$, the Cauchy problem 
\begin{equation} 
\left\{ \enspace 
\begin{aligned} 
  L U(t,x) &= F(t,x), \quad (t,x)\in (0,T)\times\Rn, \\ 
  U(0,x) &= U_0(x). 
\end{aligned} 
\right. 
\end{equation}  
possesses a solution $U\in H^{s,0;\lambda}((0,T)\times\R^{n})$. This solution
$U$ is unique in $L^2$.
\end{proposition} 
 
\begin{proof} 
We have $A_{1}=0$. Let $M_{0}$ be a symmetrizer for $A_{0}$, and 
$M_{1}=0$. Then~\eqref{PRE} is satisfied with $\delta=0$. The assertion 
follows immediately from Theorem~\ref{main}. 
\end{proof} 

 
\subsection{Characteristic roots of constant multiplicity} 
 
 
An interesting class of systems to which Theorem~\ref{main} applies is that of
systems having characteristic roots of constant multiplicity.
 
\begin{definition} 
System~\eqref{syst} is said to have {\em characteristic roots of constant
  multiplicity\/} if it is sym\-me\-triz\-a\-ble-hyperbolic in the sense of
Definition~\ref{SymmHypDef} and if 
\[
  \det \left(\tau\bsen-\sigma^{1}(A)(t,x,\xi)\right) = \prod_{h=1}^r 
  (\tau-t^\ls\mu_h(t,x,\xi))^{N_h},
\] 
where $\mu_h\in C^{\infty}([0,T],S^{(1)})$ for $1\leq h\leq r$ are
 real-valued,
 $N_{1}+\dots+N_{r}=N$, and
\[ 
|\mu_{h}(t,x,\xi)-\mu_{h'}(t,x,\xi)| 
\geq 
c, 
\quad 
1\leq h<h'\leq r,
\]
for some $c>0$. 
\end{definition} 

\begin{remark}\label{strict_hyperb}
In case $r=N$, we have $N_{1}=\dots=N_{r}=1$ and the operator
$D_{t}-A(t,x,D_x)$ is strictly hyperbolic for $t>0$.
\end{remark} 
 
If $D_{t}-A(t,x,D_{x})$ has characteristic roots of constant multiplicity,
then there exists a matrix $M_{0}\in C^{\infty}([0,T],S^{(0)})$ with $|\det
M_{0}(t,x,\xi)|\geq c>0$ such that
\[ 
B_{0}(t,x,\xi) 
:= 
(M_{0}A_{0}M_{0}^{-1})(t,x,\xi) 
= 
\diag(\mu_{1}\bsenn{N_{1}},\dots,\mu_{r}\bsenn{N_{r}})(t,x,\xi) 
\] 
is a diagonal matrix. With $A_{1}(x,\xi)=il_{\ast}^{-1}\tilde{\sigma}^{0,1}
(A)(x,\xi)$ as before, we put
\begin{equation}\label{aqw}
B_{1}(x,\xi) 
:= 
M_{0}(0,x,\xi)A_{1}(x,\xi)M_{0}(0,x,\xi)^{-1} 
= 
\begin{pmatrix} 
B_{1,11} & B_{1,12} & \dots  & B_{1,1r} \\ 
B_{1,21} & B_{1,22} & \dots  & B_{1,2r} \\   
\vdots   & \vdots   & \ddots & \vdots \\ 
B_{1,r1} & B_{1,r2} & \dots  & B_{1,rr} 
\end{pmatrix} 
\end{equation}
where $B_{1,jk}\in C^{\infty}([0,T],S^{(0)})$ is an $N_{j}\times N_{k}$
matrix.
 
\begin{proposition} 
Assume system~\eqref{syst} has characteristic roots of constant multiplicity,
and define $B_{0}$, $B_{1}$ as above. Let $\delta\in\cB(\R^{n};\R)$ be so that
\begin{equation}\label{ner}
\Re B_{1,jj}(x,\xi)\leq \delta(x)\bsenn{N_{j}}, 
\quad 
(x,\xi)\in \Rn\times(\Rn\setminus0), \quad 1\leq j\leq r. 
\end{equation}
Then, for all $s\geq 0$, $U_{0}\in \hn{s+\beta_{\ast}\delta(x)l_{\ast}}$, and
$F\in H^{s,\delta(x);\lambda}$, the Cauchy problem~\eqref{syst} possesses a
unique solution $U\in H^{s,\delta(x);\lambda}$.
\end{proposition} 
\begin{proof} 
Assuming \eqref{ner}, we are looking for a matrix $M_{1}\in S^{(0)}$ such that 
\[ 
\Re(B_{1}+\comm{M_{1}M_{0}^{-1}}{B_{0}})(0,x,\xi) 
\leq 
\delta(x)\bsen, 
\quad 
(x,\xi)\in \Rn\times(\Rn\setminus0). 
\] 
We are done if we can find a matrix $P_{1}=M_{1}M_{0}^{-1}$ in such a way that 
\begin{equation}\label{euw}
B_{1}+\comm{P_{1}}{B_{0}} 
= 
\begin{pmatrix} 
B_{1,11} & 0        & \dots  & 0        \\ 
0        & B_{1,22} & \dots  & 0        \\   
\vdots   & \vdots   & \ddots & \vdots  \\ 
0        & 0        & \dots  & B_{1,rr} 
\end{pmatrix} 
\end{equation}
is block-diagonal. 

Such a matrix $P_{1}$ can be constructed using the fact that $B_{0}$ is
diagonal with distinct eigenvalues for the different blocks, and employing the
following result, see {\sc Taylor}~\cite[Chap.~IX, Lemma~1.1]{Tay81}:

\emph{For $E\in M_{M\times M}(\C)$, $F\in M_{N\times N}(\C)$, the map
\[
  \textup{$M_{M\times N}(\C)\to M_{M\times N}(\C)$,} \quad 
  \textup{$T\mapsto TF-ET$,}
\] 
is bijective if and only if $E$ and $F$ have disjoint spectra.} 

We choose $P_1$ so that $P_{1,jj}=0$ for $1\leq j\leq r$, where the meaning of
$P_{1,jk}$ is the same as in \eqref{aqw}. Then $\comm{P_1}{B_0}_{jk} =
P_{1,jk}B_{0,kk} - B_{0,jj}P_{1,jk}$ for $j\neq k$, while
$\comm{P_1}{B_0}_{jj} = 0$ for $1\leq j\leq r$. According to the result
just quoted, we can choose $P_{1,jk}$ for $j\neq k$ so that
\[
  B_{1,jk} + \comm{P_1}{B_0}_{jk} = 0, \quad j\neq k.
\]
That is, by this choice of $P_1$ we kill all off-diagonal entries of $B_1$,
while the diagonal entries of $B_1$ remain unchanged. Thus, we end up with
\eqref{euw}.
\end{proof} 
 
\begin{example}\label{gustr}
There are two extreme cases exemplified by (a), (b) below:

(a) Let $r=N$, see Remark~\ref{strict_hyperb}. Then we can choose any
$\delta\in\cB^\infty(\Rn;\R)$ satisfying
\[ 
\Re(M_{0}A_{1}M_{0}^{-1})_{jj}(0,x,\xi)\leq \delta(x), 
\quad 
(x,\xi)\in \Rn\times(\Rn\setminus0), \quad 1\leq j\leq N. 
\] 
In a forthcoming paper, we will show that this bound on $\delta$ is sharp. 

(b) Consider the Cauchy problem 
\begin{equation} 
\left\{ \enspace 
\begin{aligned} 
  & D_t U(t,x) + i\ls a(t,x)\,h(t,D_x)U(t,x) = F(t,x),  
  \quad (t,x)\in (0,T)\times\Rn, \\ 
  & U(0,x) = U_0(x), 
\end{aligned} 
\right. 
\end{equation}  
where $a\in \cB^\infty([0,T]\times\Rn;M_{N\times N}(\C))$. Then
$A(t,x,\xi)=-il_{\ast}a(t,x)h(t,\xi)$, $A_{0}(t,x,\xi)=0$, and
$A_{1}(t,x,\xi)=a(t,x)$. By choosing $M_{0}(t,x,\xi)$ so that $M_{0}(0,x,\xi)$
is unitary and diagonalizes $\Re a(0,x)$, we see that we can choose any
$\delta\in\cB^\infty(\Rn;\R)$ satisfying
\[ 
  \delta(x) \geq \max_{1\leq j\leq N} \nu_j(x), 
\] 
where $\nu_1(x),\dots,\nu_N(x)$ are the eigenvalues of $\Re a(0,x)$ (not
necessarily distinct).
\end{example} 
 
 
\subsection{Higher-order scalar equations} 
 
 
Let $L$ be the operator 
\[ 
  L = D_t^m + \sum_{\substack{j+|\alpha|\leq m,\\j<m}} a_{j\alpha}(t,x)\, 
  t^{(j+(\ls+1)|\alpha|-m)^+}D_t^j D_x^\alpha, \quad (t,x)\in (0,T)\times\Rn, 
\] 
where $a_{j\alpha}\in \cB^\infty([0,T]\times\Rn)$ for $j+|\alpha|\leq m$, 
$j<m$. We assume $L$ to be strictly hyperbolic in the sense that 
\[ 
  \sigma^m(L) = \prod_{h=1}^m (\tau-t^\ls\mu_h(t,x,\xi)), 
\] 
where $\mu_h\in C^{\infty}([0,T],S^{(1)})$, $1\leq h\leq m$, are real--valued, 
and 
\[ 
  \bigl|\mu_h(t,x,\xi)-\mu_{h'}(t,x,\xi)\bigr| \geq c\,|\xi|, \quad  
  1\leq h <h'\leq m, 
\quad c>0. 
\] 
 
We define a reduced principal symbol of $L$, 
\[ 
  p(\tau):= p(t,x,\tau,\xi)= 
  \tau^m + p_{m-1}\tau^{m-1} + \dots + p_1\tau  + p_0, 
\] 
where 
\[ 
  p_j(t,x,\xi) := \sum_{|\alpha|=m-j}a_{j\alpha}(t,x) 
  \bigl(\frac{\xi}{|\xi|}\bigr)^\alpha, 
\] 
and a reduced secondary symbol, 
\[ 
  q(\tau):= q(x,\tau,\xi) = 
  q_{m-2}\tau^{m-2} + q_{m-3}\tau^{m-3} + \dots + q_1\tau  + q_0, 
\] 
where 
\[ 
  q_j(x,\xi) := i\ls^{-1}\sum_{|\alpha|=m-j-1}a_{j\alpha}(0,x) 
  \bigl(\frac{\xi}{|\xi|}\bigr)^\alpha. 
\] 
 
The loss of regularity is then determined as follows: 
\begin{proposition}\label{prop1} 
Let $s\geq0$, $\delta\in \cB^\infty(\Rn;\R)$ satisfy 
\begin{equation} 
  \delta(x) \geq \sup_{1\leq h\leq m}\sup_{|\xi|=1}\left(-\,\frac{ 
  \frac{\tau}{2}\,\frac{\partial^2 p}{\partial\tau^2}+ 
  \Re q}{\frac{\partial p}{\partial\tau}} 
  \right)(0,x,\mu_h(0,x,\xi),\xi). 
\end{equation} 
Then, for all $u_j\in H^{s+m-j\bs-1+\bs\delta(x)\ls}$ for $0\leq 
j\leq m-1$, $f\in H^{s,\delta(x)+m-1;\lambda}$, the Cauchy problem 
\begin{equation}\label{higher} 
\left\{\quad 
\begin{aligned} 
   Lu(t,x) & = f(t,x), \quad  (t,x)\in (0,T)\times\Rn, \\ 
   D_t^ju(0,x) &= u_j(x), \quad 0\leq j\leq m-1, 
\end{aligned} 
\right. 
\end{equation} 
possesses a solution $u\in H^{s+m-1,\delta(x);\lambda}$. This solution $u$ is 
unique in the space $H^{m-1,\delta(x);\lambda}$.  
\end{proposition} 
\begin{proof} 
We convert problem \eqref{higher} into an $m\times m$ system of the first 
order. Then it is equivalent to the Cauchy problem 
\[ 
\left\{\quad 
\begin{aligned} 
   D_t U(t,x) & = A(t,x,D_x) U(t,x) + F(t,x),  
    \quad (t,x)\in (0,T)\times\R^{n},\\ 
 U(0,x) & = U_0(x), 
\end{aligned} 
\right. 
\] 
where $U=\left(\begin{smallmatrix} g^{m-1}u \\ 
    g^{m-2}D_tu \\ \vdots \\ g D_t^{m-2}u \\ D_t^{m-1}u 
  \end{smallmatrix}\right)\in H^{s,\delta(x)+m-1;\lambda}$  
(if and only if $u\in H^{s+m-1,\delta(x);\lambda}$, see
Remark~\ref{later_use}~(b) (i)),
$U_0=\left(\begin{smallmatrix}\langle D_x\rangle^{\bs(m-1)} u_0 \\
    \langle D_x\rangle^{\bs(m-2)} u_1 \\ \vdots \\
    \langle D_x\rangle^\bs u_{m-2} \\ u_{m-1}
  \end{smallmatrix}\right)\in H^{s+\bs(\delta(x)+m-1)\ls}$, 
$F=\left(\begin{smallmatrix} 0 \\ 0 \\ \vdots \\  
    0 \\ f(t,x) 
\end{smallmatrix}\right)\in H^{s,\delta(x)+m-1;\lambda}$,  
\[ 
  A(t,x,\xi) = \begin{pmatrix} 
  (m-1)\,\frac{D_t g}{g} & g & 0 & \hdots & 0 & 0 \\ 
  0 & (m-2)\,\frac{D_t g}{g} & g & \hdots & 0 & 0 \\ 
  0 & 0 & (m-3)\,\frac{D_t g}{g} & \hdots & 0 & 0 \\ 
  \vdots & \vdots & \vdots & \ddots & \vdots & \vdots \\ 
  0 & 0 & 0 & \hdots & \frac{D_t g}{g} & g \\ 
  -\frac{a_0}{g^{m-1}} & -\frac{a_1}{g^{m-2}} & -\frac{a_2}{g^{m-3}} &  
  \hdots & -\frac{a_{m-2}}{g} & -a_{m-1}  
  \end{pmatrix}, 
\] 
and $a_j(t,x,\xi) = \sum_{|\alpha|\leq m-j} a_{j\alpha}(t,x)\,t^{(j+(\ls+1) 
  |\alpha|-m)^+}\xi^\alpha$. 
 
From Example~\ref{exam2} (b), we infer that $A\in \tilde S^{1,1;\lambda}$, 
$\sigma^1(A)(t,x,\xi)=\lambda(t)|\xi|A_0(t,x,\xi)$, where 
\[ 
  A_0(t,x,\xi) = \begin{pmatrix} 
  0 & 1 & 0 & \hdots & 0 & 0 \\ 
  0 & 0 & 1 & \hdots & 0 & 0 \\ 
  0 & 0 & 0 & \hdots & 0 & 0 \\ 
  \vdots & \vdots & \vdots & \ddots & \vdots & \vdots \\ 
  0 & 0 & 0 & \hdots & 0 & 1 \\ 
  -p_0 & -p_1 & -p_2 & \hdots & -p_{m-2} & -p_{m-1} 
  \end{pmatrix}, 
\] 
and $\tilde \sigma^{0,1}(A)(x,\xi)=-i\ls A_1(x,\xi)$, where 
\[ 
  A_1(x,\xi) =  \begin{pmatrix} m-1 & 0 & 0 & \hdots & 0 & 0 \\ 
  0 & m-2 & 0 & \hdots & 0 & 0 \\ 
  0 & 0 & m-3 & \hdots & 0 & 0 \\ 
  \vdots & \vdots & \vdots & \ddots & \vdots & \vdots \\ 
  0 & 0 & 0 & \hdots & 1 & 0 \\ 
  -q_0 & -q_1 & -q_2 & 
  \hdots & -q_{m-2} & 0 
  \end{pmatrix}. 
\] 
 
Now, it is easy to provide a symmetrizer $M_{0}$ for $A_0$, namely 
\[ 
  M_{0}(t,x,\xi)^{-1} = \begin{pmatrix} 
  1 & 1 & \hdots & 1 \\ 
  \mu_1 & \mu_2 & \hdots & \mu_m \\ 
  \vdots & \vdots & \ddots & \vdots \\ 
  \mu_1^{m-1} & \mu_2^{m-1} & \hdots & \mu_m^{m-1} 
  \end{pmatrix}. 
\] 
Note that $\det M_{0}^{-1}=\prod_{h>h'} (\mu_h-\mu_{h'})$ and, for $1\leq 
h,\,j\leq m$, 
\begin{equation}\label{help} 
  (M_{0}(t,x,\xi))_{hj} = \frac{\mu_h^{m-j} + p_{m-1}\mu_h^{m-j-1} + \dots + 
  p_{j+1}\mu_h + p_j}{\frac{\partial p}{\partial\tau}(\mu_h)}. 
\end{equation} 
 
According to our general scheme, see Example~\ref{gustr} (a), to read off the
loss of regularity we have to calculate
\begin{align*} 
  \left(M_{0}A_1M_{0}^{-1}\right)_{hh} 
  & =\sum_{j,\,k}(M_{0})_{hj}\,(A_1)_{jk}\,(M_{0}^{-1})_{kh} \\ 
  & = \sum_{j=1}^{m-1} (m-j)\, (M_{0})_{hj}\,(M_{0}^{-1})_{jh} -  
  \sum_{j=1}^{m-1} q_{j-1}\,(M_{0})_{hm}\,(M_{0}^{-1})_{jh} \\ 
  & = m - \sum_{j=1}^{m} j\, (M_{0})_{hj}\,(M_{0}^{-1})_{jh} -  
  \sum_{j=1}^{m-1} q_{j-1}\,(M_{0})_{hm}\,(M_{0}^{-1})_{jh}.  
\end{align*} 
By virtue of \eqref{help}, 
\begin{align*} 
  \sum_{j=1}^m j (M_{0})_{hj} (M_{0}^{-1})_{jh} &= \frac{1}{\frac{\partial p}{ 
  \partial\tau}(\mu_h)}\sum_{j=1}^m j\left[\mu_h^{m-j} +  
  p_{m-1}\mu_h^{m-j-1} + \dots + p_{j+1}\mu_h + p_j\right]\mu_h^{j-1} \\ 
  &= \frac{\sum_{j=1}^m\binom{j+1}{2}p_j\mu_h^{j-1}}{ 
  \frac{\partial p}{\partial\tau}(\mu_h)} 
  = \left(\frac{\frac{\partial p}{\partial\tau}+ 
  \frac{\tau}{2}\,\frac{\partial^2 p}{\partial\tau^2}}{ 
  \frac{\partial p}{\partial\tau}}\right)(0,x,\mu_h,\xi) 
\end{align*} 
and 
\[ 
  \sum_{j=1}^{m-1} q_{j-1}\,(M_{0})_{hm}\,(M_{0}^{-1})_{jh} =  
  \frac{\sum_{j=1}^{m-1}q_{j-1}\mu_h^{j-1}}{ 
  \frac{\partial p}{\partial\tau}(\mu_h)} 
  = \frac{q(x,\mu_h,\xi)}{\frac{\partial p}{\partial\tau}(0,x,\mu_h,\xi)}. 
\] 
Hence, the assertion follows. 
\end{proof} 
 
\begin{remark} 
The expression  
\[ 
  \ls\sup_{x\in\Rn,\,|\xi|=1}\left(-\,\frac{ 
  \frac{\tau}{2}\,\frac{\partial^2 p}{\partial\tau^2}+ 
  \Re q}{\frac{\partial p}{\partial\tau}} 
  \right)(0,x,\mu_h(0,x,\xi),\xi)  
\] 
is the connecting coefficient $m_h^+$ from \textsc{Amano--Nakamura} 
\cite{AmanoNakamura84}. 
\end{remark} 
 
 
\renewcommand{\thesection}{\Alph{section}} 
\setcounter{section}{0} 
\section{Appendices} 
 
\subsection{A useful estimate} 
 
 
We consider a matrix pseudodifferential operator $\del_{t}-B(t,x,D_{x})$ and 
its forward fundamental solution $X(t,t')$ which is defined by the relations 
\begin{align*} 
& 
(\del_{t}-B(t,x,D_{x}))\circ X(t,t')=0, 
\quad  
0\leq t'\leq t\leq T, \\
& 
X(t',t')=I,\quad 0\leq t'\leq T. 
\end{align*} 
We suppose that this forward fundamental solution operator exists and maps 
$\S(\Rn)$ onto itself.   
Our assumptions on $B(t,x,D_{x})$ are as follows: 
\begin{itemize} 
\item  
$B(t,x,\xi)=B_{1}(t,x,\xi)+B_{r}(t,x,\xi)$ with  
$B_{1}\in\lp{\infty}((0,T),S^{1}_{1,0})$, 
$B_{r}\in\lp{\infty}((0,T),S^{0}_{\vr,\delta})$, where 
$0\leq\delta\leq\vr\leq1$, $\delta<1$, 
\item  
$B_{1}(t,x,\xi)+B_{1}^{\ast}(t,x,\xi)\leq 2q(t,x,\xi)\bsen$ for all 
$(t,x,\xi)\in[0,T]\times\R^{2n}$, where $B_{1}^{\ast}(t,x,\xi)$ denotes the 
Hermitian conjugate of the matrix $B_{1}(t,x,\xi)$, 
\item 
the real-valued scalar function $q=q(t,x,\xi)$ belongs to 
$\lp{\infty}((0,T),S^{1}_{1,0})$ and depends either only on $(t,x)$ or only 
on $(t,\xi)$, 
\item  
$p(t,x,\xi)=\int_{0}^{t}q(t',x,\xi)\,dt'\in  
\lp{\infty}((0,T),S^{0}_{1,0})$. 
\end{itemize} 
Think of $B_{1}$ as the first-order principal symbol of $B$, which is almost
skew-symmetric (up to an integrable perturbation described by $q$), and regard
$B_{r}$ as remainder term.
 
\begin{lemma}\label{UsefulLem} 
Under these assumptions, the forward fundamental solution operator can be 
extended such as acting boundedly from $\lp{2}(\Rn)$ onto itself, 
\[ 
X\in \lp{\infty}(\lap_{+},\L(\lp{2}(\Rn))), 
\quad \lap_{+}=\{(t,t')\colon 0\leq t'\leq t\leq T\}. 
\] 
\end{lemma} 
 
\begin{proof} 
For $(t,t')\in\lap_{+}$, we define a mapping  
$Y(t,t')\colon \S(\Rn)\to\S(\Rn)$ by 
\[ 
Y(t,t') 
= 
\exp(-p(t,x,D_{x}))\circ \exp(p(t',x,D_{x})) \circ X(t,t'). 
\]   
Obviously, $Y(t',t')=I$. Since the symbol $p(t,x,\xi)$ does not depend on $x$ 
and $\xi$ simultaneously, we have 
\begin{align*} 
\del_{t}\circ Y(t,t')   
& 
= 
-q(t,x,D_{x})\circ Y(t,t') 
\\ 
& 
\quad\quad\quad\quad 
+ 
\exp(-p(t,x,D_{x}))\circ \exp(p(t',x,D_{x})) \circ B(t,x,D_{x}) \circ X(t,t') 
\\ 
& 
= 
\left( 
B-q\bsen+\comm{e^{-p(t,x,D_{x})}e^{p(t',x,D_{x})}\bsen}{B}
e^{-p(t',x,D_{x})}e^{p(t,x,D_{x})}   
\right)\circ Y(t,t') 
\\ 
& 
= 
\left( 
B_{1}-q\bsen+B_{0}   
\right)\circ Y(t,t') 
\end{align*} 
for some $B_{0}\in\lp{\infty}(\lap_{+},S^{0}_{\vr,\delta})$ because of 
$\exp(\pm p(t,x,\xi))\in\lp{\infty}((0,T),S^{0}_{1,0})$. 
 
For fixed $t'\in[0,T]$, $U_{0}\in\S(\Rn)$, we define a function 
$U(t,x)=Y(t,t')U_{0}(x)$ which solves 
\begin{align*}
& \del_{t}U=(B_{1}-q\bsen +B_{0})U, \quad (t,x)\in (t',T)\times\R^{n}, 
\\ 
& U(t',x)=U_{0}(x). 
\end{align*}
Employing the sharp G\r{a}rding inequality and Calder\'on-Vaillancourt's
theorem, we obtain
\begin{multline*} 
\del_{t}\SP{U(t,\cdot)}{U(t,\cdot)} = 
2\Re\SP{\del_{t}U(t,\cdot)}{U(t,\cdot)} 
= 
2\Re\SP{(B_{1}-q\bsen+B_{0})U(t,\cdot)}{U(t,\cdot)} 
\\ 
\leq 
\SP{(B_{1}+B_{1}^{\ast}-2q\bsen)U(t,\cdot)}{U(t,\cdot)} 
+ 
2\norm{(B_{0}U)(t,\cdot)}{\lp{2}} 
\norm{U(t,\cdot)}{\lp{2}} 
\leq 
C\norm{U(t,\cdot)}{\lp{2}}^{2}.
\end{multline*} 
Then Gronwall's lemma implies  
$\norm{U(t,\cdot)}{\lp{2}}\leq C\norm{U(t',\cdot)}{\lp{2}}$, i.e.,
\[ 
Y\in\lp{\infty}(\lap_{+},\L(\lp{2}(\R^{n}))). 
\] 
The operators $\exp(\pm p(t,x,D_x))$ map $\lp{2}(\R^{n})$ continuously and
bijectively onto itself which completes the proof.
\end{proof} 

 
\subsection{Proof of Proposition~\ref{asd} (b)} 
 

We need the following result:

\begin{lemma}\label{cici}
For each $N\times N$ matrix function $q_0\in S^{(0)}$ satisfying $|\det
q_0(x,\xi)|\geq c$ for all $(x,\xi)\in\Rn\times(\Rn\setminus0)$ and some $c>0$,
there is an invertible operator $Q\in S_\cl^0(\Rn)$ such that
$\sigma^0(Q)(x,\xi)=q_0(x,\xi)$.
\end{lemma}
\begin{proof}
We construct two invertible operators $Q_1,\,Q_2\in S_\cl^0(\Rn)$ such that
\[ 
  \sigma^0(Q_1)(x,\xi) = q_0(x,\xi)q_0(x^0,\xi)^{-1}, \quad
  \sigma^0(Q_2)(x,\xi) = q_0(x^0,\xi)
\]
for some $x^0\in \Rn$. Then the composition $Q_1Q_2$ has the desired
properties.

\emph{Construction of\/ $Q_1$.} \enspace 
We employ the parameter-dependent calculus of \textsc{Grubb} \cite{Gr}.

Rename $q_0(x,\xi)q_0(x^0,\xi)^{-1}$ to $q_0(x,\xi)$. Then $q_0(x^0,\xi)=
\bsen$ for all $\xi\in\Rn\setminus0$. Therefore, there is an invertible 
$N\times N$ matrix function $p_0\in S^{(0)}\bigl(\Rn\times((\Rn\times
\overline{\R}_+)\setminus0)\bigr)$ such that $|\det p_0(x,\xi,\mu)|\geq c/2$
for $(x,\xi,\mu)\in \Rn\times ((\Rn\times \overline{\R}_+)\setminus0)$ and
\[
  p_0(x,\xi,0) = q_0(x,\xi), \quad (x,\xi)\in\Rn\times(\Rn\setminus0).
\]
We now set
\[
  p(x,\xi,\mu) := \chi(|\xi,\mu|)\left(p_1(x,\xi,\mu) + \chi(|\xi|)
  \bigl(p_0(x,\xi,\mu)-p_1(x,\xi,\mu)\bigr)\right),
\]
where $p_1(x,\xi,\mu):=\sum_{|\alpha|<k}\frac{\xi^\alpha}{\alpha!}\,
\partial_\xi^\alpha p_0(x,\xi,\mu)$ for some integer $k>0$, see \cite[Remark
2.1.13]{Gr}. According to \cite[Theorem 3.2.11]{Gr}, there is a $\mu_0\geq0$
such that, for all $\mu\geq\mu_0$, the operator $p(x,D_x,\mu)\colon
L^2(\Rn)\to L^2(\Rn)$ is invertible. It suffices to set $Q_1:=p(x,D_x,\mu)$,
where $\mu\geq\mu_0$.

\emph{Construction of\/ $Q_2$.}\enspace Rename $q_0(x^0,\xi)$ to $q_0(\xi)$.
The task to construct $q\in S_\cl^0$ such that $\sigma^0(q)(x,\xi)=
q_0(\xi)$ and $q(x,D_x)\in \Op S_\cl^0$ is invertible can be fulfilled within
the framework of $SG$-calculus, where one has symbols which have asymptotic
expansions into components which are homogeneous in both the $x$- and the
$\xi$-variables. In particular, we have a symbol $\sigma_e^0(q)(x,\xi)\in
S^{(0)}(\Rn_x\setminus0) \hat\otimes_\pi S_\cl^0(\Rn_\xi)$, having the
status of a second principal symbol, subject only to the restriction $\sigma^0
(\sigma_e^0(q)(x,\xi))=q_0(\xi)$. Choosing $\sigma_e^0(q)(x,\xi)$ as an
elliptic symbol in $x\neq0$ uniformly with respect to $\xi$, we get that
$q(x,D_x)\colon L^2(\Rn)\to L^2(\Rn)$ is a Fredholm operator. Moreover, upon
an appropriate choice of $\sigma_e^0(q)(x,\xi)$ we can achieve each integer as
index of this operator. We choose $\sigma_e^0(q)(x,\xi)$ in such a way that
$q(x,D_x)\colon L^2(\Rn)\to L^2(\Rn)$ has index $0$. Then, by adding a
contribution from $\Op S^{-\infty}(\Rn_x\times\Rn_\xi)$ if necessary, we
finally arrive at an operator $q(x,D_x)$ that is invertible as operator on
$L^2(\Rn)$. We leave the details of this construction to the reader. For more
on $SG$-calculus we refer, e.g., to \textsc{Schulze} \cite{Schulze98}.
\end{proof}

\begin{proof}[Proof of\/ \textup{Proposition~~\ref{asd} (b)}]
There is a generalization of Lemma~\ref{cici} to the case $q_0\in
C^\infty([0,T];S^{(0)})$. Therefore, we find an invertible operator $Q_1\in
C^\infty([0,T];\Op S^{(0)})$ such that $\sigma^0(Q_1)(t,x,\xi)=q_0(t,x,\xi)$
for $(t,x,\xi)\in [0,T]\times\Rn\times(\Rn\setminus0)$. Since $q\in
C^\infty([0,T];S_\cl^0+S^{-1})$ implies $q\in \tilde S^{0,0;\lambda}$, where
$\tilde\sigma^{-1,0}(q)=0$, it remains to construct an operator $Q_2\in
\Op\tilde S^{0,0;\lambda}$ such that
\[
  \sigma^0(Q_2)(t,x,\xi) = \bsen, \quad 
  \tilde\sigma^{-1,0}(Q_2)(x,\xi) =q_0(0,x,\xi)^{-1}q_1(x,\xi),
\]
and the composition $Q_1Q_2$ has the desired properties.

Rename $(q_0^{-1}q_1)(0,x,\xi)$ to $q_1(x,\xi)$ and set $Q_2=q(t,x,D_x)$,
where 
\[
  q(t,x,\xi) = \bsen + \chi(\Lambda(t)\lxr/d) t^{-(\ls+1)} q_1(x,\xi)
\] 
for some large $d>0$ to be chosen. We have
\[
  \left|\chi(\Lambda(t)\lxr/d) t^{-(\ls+1)} q_1(x,\xi)\right| \leq 
  Cd^{-1}, \quad (t,x,\xi)\in [0,T]\times\R^{2n},
\]
for some $C>0$ and $d>0$ is large enough. From \textsc{H\"ormander}
\cite[Theorem 18.1.15]{Hoe85}, we conclude that
\[
  \left\|\chi(\Lambda(t)\lDr/d) t^{-(\ls+1)} q_1(x,D_x)\right\|_{
  \cL(L^2(\Rn))} \leq \frac13 + C'\left(\frac{\Lambda(t)}{d}\right)^{1/2}, 
  \quad t\in[0,T],
\]
for some $C'>0$ and $d\geq C/3$ is large enough. Choosing additionally $d\geq
9C'\,{}^2\Lambda(T)$, we find that, for each $t\in[0,T]$, the operator
$q(t,x,D_x)$ is invertible on $L^2(\Rn)$ with
\[
  \|q(t,x,D_x)^{-1}\|_{\cL(L^2(\Rn))}\leq 3.
\] 
This completes the proof.
\end{proof}
 

 
\bibliographystyle{amsplain} 
\bibliography{energy_estimates} 

\providecommand{\bysame}{\leavevmode\hbox to3em{\hrulefill}\thinspace}
\providecommand{\MR}{\relax\ifhmode\unskip\space\fi MR }
\providecommand{\MRhref}[2]{%
  \href{http://www.ams.org/mathscinet-getitem?mr=#1}{#2}
}
\providecommand{\href}[2]{#2}
\begin{thebibliography}{10}

\bibitem{AmanoNakamura84}
K.~Amano and G.~Nakamura, \emph{Branching of singularities for degenerate
  hyperbolic operators}, Publ. Res. Inst. Math. Sci. \textbf{20} (1984),
  225--275.

\bibitem{BdM74}
L.~Boutet~de Monvel, \emph{Hypoelliptic operators with double characteristics
  and related pseudo-differential operators}, Comm. Pure Appl. Math.
  \textbf{27} (1974), 585--639.

\bibitem{DreherOsaka}
M.~Dreher, \emph{Weakly hyperbolic equations, {S}obolev spaces of variable
  order, and propagation of singularities}, Osaka J. Math. \textbf{39} (2002),
  409--445.

\bibitem{DreherReissig}
M.~Dreher and M.~Reissig, \emph{Propagation of mild singularities for
  semilinear weakly hyperbolic equations}, J. Analyse Math. \textbf{82} (2000),
  233--266.

\bibitem{DW04}
M.~Dreher and I.~Witt, \emph{Parametrix construction for weakly hyperbolic
  operators}, In preparation.

\bibitem{DW02}
\bysame, \emph{Edge {S}obolev spaces and weakly hyperbolic equations}, Ann.
  Mat. Pura Appl. \textbf{180} (2002), 451--482.

\bibitem{Gr}
G.~Grubb, \emph{Functional calculus of pseudo-differential boundary problems},
  second ed., Progr. Math., vol.~65, Birkh\"auser Boston, Boston, MA, 1996.

\bibitem{Han79}
N.~Hanges, \emph{Parametrices and propagation for operators with non-involutive
  characteristics}, Indiana Univ. Math. J. \textbf{28} (1979), 87--97.

\bibitem{Hoe85}
L.~H\"ormander, \emph{The analysis of linear differential operators {III}},
  Grundlehren Math. Wiss., vol. 274, Springer, Berlin, 1985.

\bibitem{IP74}
V.~Ivrii and V.~Petkov, \emph{Necessary conditions for the {C}auchy problem for
  non-strictly hyperbolic equations to be well-posed}, Russian Math. Surveys
  \textbf{29} (1974), 1--70.

\bibitem{Jos98}
M.~S. Joshi, \emph{A symbolic construction of the forward fundamental solution
  of the wave operator}, Comm. Partial Differential Equations \textbf{23}
  (1998), 1349--1417.

\bibitem{NakamuraUryu}
G.~Nakamura and H.~Uryu, \emph{Parametrix of certain weakly hyperbolic
  operators}, Comm. Partial Differential Equations \textbf{5} (1980), 837--896.

\bibitem{Qi}
M.-Y. Qi, \emph{On the {C}auchy problem for a class of hyperbolic equations
  with initial data on the parabolic degenerating line}, Acta Math. Sinica
  \textbf{8} (1958), 521--529.

\bibitem{Schulze98}
B.-W. Schulze, \emph{Boundary value problems and singular pseudo-differential
  operators}, Wiley Ser. Pure Appl. Math., J. Wiley, Chichester, 1998.

\bibitem{TaniguchiTozaki}
K.~Taniguchi and Y.~Tozaki, \emph{A hyperbolic equation with double
  characteristics which has a solution with branching singularities}, Math.
  Japon. \textbf{25} (1980), 279--300.

\bibitem{Tay81}
M.~Taylor, \emph{Pseudodifferential operators}, Princeton Math. Ser., vol.~34,
  Princeton Univ. Press, Princeton, NJ, 1981.

\bibitem{Wit03}
I.~Witt, \emph{A calculus for a class of finitely degenerate pseudodifferential
  operators}, Evolution Equations: Propagation Phenomena -- Global Existence --
  Influence of Non-linearities (R.~Picard, M.~Reissig, and W.~Zajaczkowski,
  eds.), Banach Center Publ., vol.~60, Polish Acad. Sci., Warszawa, 2003,
  pp.~161--189.

\bibitem{YagdjianBuch}
K.~Yagdjian, \emph{The {C}auchy problem for hyperbolic operators.}, Math.
  Topics, vol.~12, Akademie Verlag, Berlin, 1997.

\bibitem{Yos77}
A.~Yoshikawa, \emph{Construction of a parametrix for the {C}auchy problem of
  some weakly hyperbolic equation {I}}, Hokkaido Math. J. \textbf{6} (1977),
  313--344, \emph{II.} Hokkaido Math. J. \textbf{7} (1978), 1--26. \emph{III.}
  Hokkaido Math. J. \textbf{7} (1978), 127--141.

\end{thebibliography}
 
 
\end{document}